\documentclass[a4paper, 11pt]{article} 

\usepackage{amsfonts}
\usepackage{amsmath}
\usepackage{enumerate}
\usepackage[frenchb]{babel}
\usepackage[applemac]{inputenc}
\usepackage{graphicx}

\newtheorem{lem}{Lemma}
\newtheorem{prop}{Proposition}
\newtheorem{corol}[prop]{Corollary}
\newtheorem{thm}{Theorem}
\newenvironment{dem}{\noindent\textbf{Proof.}}{\hfill \framebox[0.6em]{}}
\newenvironment{rmq}{\smallskip \noindent\textbf{Remark.}}{}

\newcommand{\N}{\mathbb{N}}
\renewcommand{\P}{\mathbb{P}}
\newcommand{\R}{\mathbb{R}}
\newcommand{\E}{\mathbb{E}}

\newcommand{\1}{\mathbf{1}}

\usepackage[bf]{caption}

\title{Iterated stochastic processes : simulation and relationship with high order partial differential equations.}
\author{Mich\`ele Thieullen\thanks{Laboratoire de Probabilit\'es et Mod\`eles Al\'eatoires, Universit\'e Pierre et Marie Curie,  4, Place Jussieu, 75252 Paris cedex 05, France}{ } \thanks{michele.thieullen@upmc.fr}\\
\and Alexis Vigot $\text{}^\star$ \thanks{alexis.vigot@courriel.upmc.fr }\\
}

  \date{\today}

\begin{document}

\maketitle
\bigskip
{\small In this paper, we consider the composition of two independent processes: one process corresponds to position and the other one to time. Such processes will be called {\it iterated processes}. We first propose an algorithm based on the Euler scheme to simulate the trajectories of the corresponding iterated processes on a fixed time interval. 
This algorithm is natural and can be implemented easily. We show that it converges almost surely, uniformly in time, with a rate of convergence of order $1/4$ and propose an estimation of the error. 
We then extend the well known Feynman-Kac formula which gives a probabilistic representation of partial differential equations (PDEs), to its higher order version using iterated processes. In particular we consider general position processes which are not necessarily Markovian or are indexed by the real line but real valued. We also weaken some assumptions from previous works. We show that intertwining diffusions are related to transformations of high order PDEs. Combining our numerical scheme with the Feynman-Kac formula, we simulate functionals of the trajectories and solutions to fourth order PDEs that are naturally associated to a general class of iterated processes.

\bigskip
\noindent \underline{Keywords} : Iterated process, Euler scheme, high order partial differential equation, Feynman-Kac formula, diffusion processes, iterated Brownian motion.}

\newpage

\section{Introduction}
In the present paper we address {\it iterated processes} for which we propose a numerical scheme to simulate their trajectories and investigate their relationship with high order partial differential equations (PDEs). Then we use our scheme to approximate numerically the solution of such PDEs. Consider two independent processes $(X_t)_{t\geq 0}$ and $(Y_t)_{t\in I}$, where $I$ is an interval, defined on a probability space $(\Omega, \mathcal{F}, \P)$ and taking values respectively in $\R^d$ and $\R^n$. We define an interated process $Z$ by replacing the time index of $X$ by $|Y|$ where $|\cdot|$ denotes a norm on $\R^n$ as follows
\begin{equation}
\label{processusitere}
Z_t:=X(|Y_t|),\quad t\in I.
\end{equation}
In the sequel we call $(X_t)$ (resp. $(Y_t)$) the {\it position} (resp. {\it time}) process. Different types of iterated processes have been considered recently for instance by Allouba \cite{allouba1, allouba}, Burdzy \cite{burdzy1, burdzy3, burdzy2}, Khoshnevisan and Lewis \cite{khoshnevisan}. 
Some of these authors extend $X$ to a process indexed by the real line as follows
\begin{equation}
\label{realline}
X_2(t):=\left \lbrace 
\begin{array}{lcl} 
X(t)&\text{if}& t\geq 0\\ 
X'(-t)&\text{if}& t<0\\ 
\end{array}\right.
\end{equation}
where $X'$ is an independent copy of $X$. The corresponding iterated process is then defined by $X_2(Y_t)$. 

\smallskip

\noindent In order to study properties of iterated processes, in particular to be able to simulate solutions of related PDEs, we would like to simulate their trajectories and to control the error.  We haven't found any reference on the simulation of such processes in the literature. In the first part of the paper we propose a scheme based on the classical Euler scheme to simulate the trajectories of the iterated process $(Z_t)$ when $X$ and $Y$ are two diffusions. For every $T>0$ and positive integer $n$, we construct a process $\tilde{Z}^n$ piecewise linear on a regular partition of $[0,T]$ with mesh $T/n$. We will see that this scheme converges a.s. uniformly as $n$ tends to infinity on $[0,T]$ to $(Z_t)$ with a rate of convergence of order $1/4$, which means that
\begin{equation}
\label{rate}
\forall\ 0<\alpha<1/4,\quad \lim_{n\rightarrow +\infty}  \, \, \, n^\alpha\sup_{t\in[0,T]} |\tilde{Z}^n(t)-Z(t)|=0,\, \, \, \quad {\rm a.s.}
\end{equation}
It seems that the order of convergence in (\ref{rate}) for iterated processes satisfying our assumptions cannot be better than $1/4$. Indeed for a Brownian motion $Y$, we have for all $\alpha>1/2$ (see \cite{faure})
\begin{equation*}
\limsup_n  \, \, \, n^\alpha \, |Y_T-Y_T^n|=\infty, \, \, \, \, \, {\rm a.s.}
\end{equation*}
The convergence of the scheme being uniform, we can use it to simulate quantities depending on the whole trajectory of the process, not only on its value at a fixed time, like for instance the variations, the mean of a functional of the process or the measure of a subset of the function space $[0,T]^\R$ with respect to the law of the process. Our scheme is easy to implement since it only requires the simulation of $2n$ Gaussian random variables and so it is adapted to methods using simulations of many trajectories such as Monte Carlo methods. It can be extended to processes given by construction (\ref{realline}). We have chosen the Euler scheme because of its simplicity. The Milshtein approximation for instance, doesn't provide a better rate of convergence and is more complex.

\smallskip

\noindent In the subsequent part of the paper we address the connection between iterated processes and high order PDEs. There are in the literature several papers which prove that the density of specific iterated processes satisfies a Fokker-Planck type PDE (\cite{orsingher2}, \cite{orsingher}), establish a Feynman-Kac type formula (\cite{allouba}, \cite{funaki}, \cite{erkan}) for such processes or associate them to fractional equations (\cite{baeumer}, \cite{erkan2}, \cite{orsingher3}). We extend these works to general position processes in particular not necessarily Markovian (cf. section \ref{gyongy}).\\
When the time process is a Brownian motion or an $\alpha$-stable process, the equation we obtain reduces to the equation obtained respectively by \cite{allouba} and \cite{erkan}. However our result holds under more general assumptions on the initial value. Moreover in this particular case we prove that the solution of the PDE is unique. It is difficult to extend the results of this paper in a general setting when the density of the time process satisfies a PDE whose coefficients depend on the spatial variable. We provide examples of such time processes. \\
We also consider transformations going from one high order PDE to another one and we point out the relationship with intertwining diffusions. However, these results (as well as those in \cite{allouba} and \cite{erkan}) concern PDEs which contain the initial value and some of its derivatives. The PDE obtained in \cite{funaki} does not have this drawback but the underlying iterated process takes values in the complex plane. Extending the construction (\ref{realline}), we obtain a Feynman-Kac formula where the initial value does not appear any longer in the PDE. Moreover the underlying process is real valued. An application of our result provides a stochastic representation of a solution of the Euler-Bernouilli beam equation when one considers iteration of a Brownian motion by an independent Cauchy process.

\smallskip

\noindent In the last part of the paper, we implement the algorithm for the iterated Brownian motion (IBM) given by (\ref{processusitere}) when $X$ and $Y$ are two independent Brownian motions. We simulate its third and fourth order variations and the solution of the corresponding fourth order PDE.

\bigskip

\noindent The present paper is organized as follows. In section \ref{algorithme} we describe our scheme, we prove its uniform convergence and study the error. Section \ref{pdes} is devoted to the relationship between iterated processes and PDEs in the spirit of \cite{allouba} and \cite{erkan}. Transformations of PDEs are addressed in section \ref{transformations}. In section \ref{funakipde}, inspired by \cite{funaki}, we work with construction (\ref{realline}) in order to obtain PDEs that do not contain the initial value. Section \ref{application} implements our scheme in the IBM case. The Appendix (section \ref{appendix}) contains some auxiliary proofs.

\section{Numerical scheme.}\label{algorithme}

\noindent In this section we describe our numerical scheme and study its convergence. We consider two independent time inhomogeneous diffusion processes $X$ and $Y$. Our scheme is based on the idea that a pathwise approximation of $X(|Y|)$ may be obtained by the composition of the respective Euler approximations of $X$ and $Y$. More precisely for $T>0$ and a positive integer $n$, the Euler scheme for $Y$ on $[0,T]$ yields a continuous approximation $Y^n$. If $M_n$ denotes the maximum of $|Y^n|$ on $[0,T]$, the Euler scheme for $X$ on $[0,M_n]$ provides an approximation $X^n$. We prove below that a piecewise linear approximation $\tilde{Z}^n$ of the composition $Z^n:=X^n(|Y^n|)$ converges a.s. uniformly on $[0,T]$ as $n$ tends to infinity to $Z=X(|Y|)$ with rate of order $\frac{1}{4}$. 

\smallskip

\noindent We start with the stochastic differential equations (SDEs) satisfied by the position and time processes. Let us assume that the position process $(X_t)_{t\geq 0}$ starting from $X_0$ satisfies 
\begin{equation}
\label{edsX}
dX_t=b_X(t,X_t)dt+\sigma_X(t,X_t)dW_t, \ t\geq 0,
\end{equation}
where $b_X: \R \times \R^d \rightarrow \R^d$, $\sigma_X: \R \times \R^d \rightarrow \R^{d\times p}$ and $(W_t)$ is a standard $\R^p$ valued Brownian motion. The time process $(Y_t)_{t\in [0,T]}$ starting from $Y_0$ solves
\begin{equation}
\label{edsY}
dY_t=b_Y(t,Y_t)dt+\sigma_Y(t,Y_t)dW'_t, \ t\in[0,T],
\end{equation}
with $b_Y: \R \times \R^n \rightarrow \R^n$, $\sigma_Y: \R \times \R^n \rightarrow \R^{n\times q}$ and $(W'_t)$ a standard $\R^q$ valued Brownian motion independent from $W$. Several assumptions will be used throughout the paper. Regarding the starting points, we assume that

\begin{equation}
\label{initial_val}
\forall p\geq 1, \, \, \E|X_0|^p+\E|Y_0|^p<\infty.
\end{equation}
The coefficients of (\ref{edsX}) and (\ref{edsY}) are supposed to enjoy Lipschitz continuity in space as well as at most linear growth and H\" older continuity in time. For simplicity, we write down these assumptions for $b_X$ and $\sigma_X$ only. For $b_Y$ and $\sigma_Y$, the only difference is that the time variable in (\ref{edsY}) belongs to the compact interval $[0,T]$. We assume that there exist two positive real numbers $K $ and $\beta$ such that
\begin{eqnarray}
&\forall t\geq 0,& \forall x,y\in \R^d,\notag\\
&&|b_X(t,x)-b_X(t,y)|+|\sigma_X(t,x)-\sigma_X(t,y)|\leq K|x-y|, \label{lipschitz_cond}\\
&&|b_X(t,x)|^2+|\sigma_X(t,x)|^2\leq K^2(1+|x|^2), \label{growth_cond}
\end{eqnarray}
\begin{eqnarray}
\label{holder_cond}
&\forall s, t\geq 0, &\forall x\in \R^n,\notag\\
&&|b_X(s,x)-b_X(t,x)|+|\sigma_X(t,x)-\sigma_X(s,x)|\leq K|t-s|^\beta.
\end{eqnarray}

\noindent Under these assumptions equations (\ref{edsX}) and (\ref{edsY}) have unique strong solutions (cf. \cite{kar}).

\bigskip

\noindent Let us recall that the Euler scheme for $Y$ on $[0,T]$ given the regular subdivision $t_0=0<t_1<\ldots<t_n=T$ with mesh $\Delta_T=T/n$ is defined recursively by $Y_0^n:=Y_0$ and
\begin{equation}
\label{eulerpourY}
\forall t\in ]t_k, t_{k+1}], Y_t^n:=Y_{t_k}^n+\sigma_Y(t_k,Y_{t_k}^n)(W'_t-W'_{t_k})+b_Y(t_k,Y_{t_k}^n)(t-t_k).
\end{equation}
Let us define $M_n:=\sup_{t\in [0,T]} |Y_t^n|$, which is a.s. finite (cf. Proposition \ref{estimeesclassiques})  and denote by $X^n$ the approximation resulting from the Euler scheme for $X$ on $[0,M_n]$. We propose to approximate $X(|Y_t|)$ by $X^n(|Y^n_t|)$. The following statement provides an estimate of the mean error associated to this approximation, uniformly on $t\in [0,T]$. In the following, $n$ is assumed to be such that $\Delta_T=T/n\leq 1$.

\begin{thm}
\label{meanerror}
Let us fix $T>0$ and assume that $X$ and $Y$ satisfy (\ref{edsX})-(\ref{holder_cond}). Let $n\geq 1$ be an integer. For $t\in [0,T]$ we set $Z_t:=X(|Y_t|)$ and $Z_t^n:=X^n(|Y_t^n|)$.
If $\sup_{t\in [0,T]} \max \{|Y_t|, |Y_t^n|\}$ has finite Laplace transform, 
\begin{equation}\label{errorestimate}
\E (\sup_{t\in [0,T]} |Z_t-Z_t^n|^{2p})\leq C(\Delta_T)^\rho,
\end{equation}
where $C$ is a constant depending only on $K$, $p$, $T$ and $\E|X_0|^{2p}$ and $\rho$ is defined by $\rho=(p-1)\inf \{1,2\beta\}/2$. In particular (\ref{errorestimate}) holds true if $Y_0^2$ has a finite Laplace transform.
\end{thm}

\smallskip

\noindent In Theorem \ref{meanerror}, we have used the exact composition of the respective Euler schemes of $X$ and $Y$ to approximate $X(|Y|)$ by $X^n(|Y^n|)$. Actually another approximation easier to implement may be chosen. Indeed, let us consider $\bar{Y}_t^n$, a step function defined on $[0,T]$ by $\bar{Y}_t^n=Y^n(kT/n)$ for $t\in[kT/n,(k+1)T/n]$ and $\bar{X}^n$, a step function defined on $[0,M_n]$ by $\bar{X}_t^n=X^n(kM_n/n)$ for $t\in[kM_n/n,(k+1)M_n/n]$. Let $\tilde{Z}^n$ be the linear interpolation between the points $\bar{X}^n(\bar{Y}^n(kT/n))$ for $k=0,1,\ldots,n$. The computation of $\tilde{Z}^n$ is facilitated by the use of piecewise constant processes whose composition is simpler than the composition of piecewise linear ones. We now show that $Z^n$ and $\tilde{Z}^n$ both converge uniformly to $Z$ with order $\frac{1}{4}$. 

\begin{thm}
\label{rateofconvergence}
Let $X$ and $Y$ satisfy (\ref{edsX})-(\ref{holder_cond}). Then $\forall \, \, 0\leq \alpha < \frac{1}{4},$
\begin{equation}
\lim_{n\rightarrow +\infty} \, \, \, n^\alpha \, sup_{t\in [0,T]} {(|Z_t^n-Z_t|+|\tilde{Z}_t^n-Z_t|)}=0, \, \, \, \, {\rm a.s.}
\end{equation}
\end{thm}

\noindent The proof of Theorem \ref{meanerror} relies on the following proposition. 
 
\begin{prop}(cf. \cite{faure})
\label{errorboundforX}
Under assumptions (\ref{edsX})-(\ref{holder_cond}),
\begin{equation}\label{explicitbound}
\E (\sup_{t\in [0,M]} |X_t-X_t^n|^{2p})\leq C(1+M^p)M(\Delta_M)^\gamma e^{CM}
\end{equation}
where $\Delta_M:=M/n$, $\gamma:=p\sup \{1,2\beta\}$ if $\Delta_M\geq1$ and $\gamma:=p\inf \{1,2\beta\}$ otherwise. $C$ is a constant depending only on $K$, $p$, $n$ and $\E|X_0|^{2p}$.
\end{prop}
In the sequel we need the explicit dependance of (\ref{explicitbound}) on the parameter $M$, since the interval on which we approximate $X$ is random. For the sake of completeness, we provide the proof of Proposition \ref{errorboundforX} in the Appendix (section 5).

\medskip

\noindent\textbf{Proof of Theorem \ref{meanerror}.}

In the proof $C$ and $\tilde C$ denote constants depending only on $K$, $p$, $n$, $T$ and $\E|X_0|^{2p}$ which may vary from line to line. Remember that $Z_t=X(|Y_t|)$ and $Z^n_t=X^n(|Y_t^n|)$. We start with
\begin{multline} 
\E\left(\sup_{t\in[0,T]} |Z^n_t-Z_t|^{2p}\right)\leq 2^{2p-1}\E\left(\sup_{t\in[0,T]} |X(|Y_t^n|)-X(|Y_t|)|^{2p}\right) \\
+2^{2p-1}\E\left(\sup_{t\in[0,T]} |X^n(|Y_t^n|)-X(|Y_t^n|)|^{2p}\right). \label{eq:82}
\end{multline}
Let us bound the first term on the RHS. Define $M:=\sup_{t\in[0,T]} \max \{|Y_t|,|Y_t^n|\}$ and $\delta:=\sup_{t\in[0,T]} |Y_t-Y_t^n|$. Let us recall that if $M$ were fixed, Garsia-Rodemich-Rumsey Lemma (cf. \cite{garsia}) would provide a random variable $\Gamma$ such that $|X_s-X_t|^{2p}\leq c\, \Gamma\,  |t-s|^m,$ for all $s,t$ in $[0,M]$ whatever $m\in]0,p-1[$. In this inequality $c$ is a constant depending on $p$ and $m$. Moreover this lemma provides the following expression for $\Gamma$,
\begin{equation}\label{expressiondegamma}
\Gamma= \int_{[0,M]^2}\, \frac{|X_a-X_b|}{|a-b|^{m+2}}\, da db,
\end{equation}
as well as the estimate
\begin{equation*}
\E (\Gamma)\leq H c_1 \frac{1}{(p-1-m)}(\frac{M}{2})^{p-m},
\end{equation*}
where $c_1$ is a universal constant and $H$ denotes the right-hand side of (\ref{Kolmogorovcriterion}). Here the pair $(M,\delta)$ is random. However by definition it is independent of $X$. Therefore we can condition on $(M,\delta)$ and write 
\begin{eqnarray*}
\E\left(\sup_{t\in[0,T]} |X(|Y_t^n|)-X(|Y_t|)|^{2p}\right)&\leq& {\tilde C}\, \E\left[ (1+M^p)e^{CM}(\frac{M}{2})^{p-m}\delta^m \right]\\
&\leq& \tilde C\left(\E\left[ P(M)e^{CM}\right] \right)^{1/2} (\E [\delta^{2m} ])^{1/2},
\end{eqnarray*}
In the latter expression (obtained by H\"older inequality) we have set $P(x):=(1+x^p)^2(\frac{x}{2})^{2(p-m)}$. Remember now that $\delta$ is equal to $\sup_{t\in[0,T]} |Y_t-Y_t^n|$. Proposition 14 in \cite{faure} implies
\begin{equation*}
\E [\delta^{2m} ]\leq (\Delta_T)^{m\inf (2\beta,1)}.
\end{equation*}
Therefore we obtain the following upper bound for the first term: 
\begin{equation}
\label{eq:80}
\E\left(\sup_{t\in[0,T]} |X(|Y_t^n|)-X(|Y_t|)|^{2p}\right)\leq \tilde C\left(\E\left[ P(M)e^{CM}\right] \right)^{1/2}(\Delta_t)^{(p-1)\inf \{2\beta,1\}/2}.
\end{equation}
The second term can be easily dominated as follows. Indeed notice first that $\sup_{t\in[0,T]} |X^n(|Y_t^n|)-X(|Y_t^n|)|\leq \sup_{t\in[0,M]} |X_t^n-X_t|$. Conditioning by $M$ (which is independent of $X$) and applying Proposition \ref{errorboundforX} we obtain,
\begin{align}
\E\left(\sup_{t\in[0,T]} |X^n(|Y_t^n|)-X(|Y_t^n|)|^{2p}\right)&\leq \E\left(\sup_{t\in[0,M]} |X_t^n-X_t|^{2p}\right)\notag \\
&\leq n^{-\gamma}\, \, \E \left[C(1+M^p)M^{1+\gamma} e^{CM}\right]. \label{eq:81}
\end{align}
where $\gamma=p\inf \{1,2\beta\}$. Combining (\ref{eq:82}) with inequalities (\ref{eq:80}) and (\ref{eq:81}) gives the result.

\noindent Let us now prove that the assumption $\sup_{t\in [0,T]} \max \{|Y_t|, |Y_t^n|\}$ has finite Laplace transform is satisfied if $Y_0^2$ has a finite Laplace transform.\\
\noindent Take $\lambda>0$. Then
$$\E e^{\lambda \sup_{t\in[0,T]} \max \{|Y_t|,|Y_t^n|\}}\leq \E e^{\lambda \sup_{t\in[0,T]} |Y_t|}+\E e^{\lambda \sup_{t\in[0,T]} |Y_t^n|},$$
\begin{align*}
\E e^{\lambda \sup_{t\in[0,T]} |Y_t|}&\leq e^\lambda +\E \left( \sum_{k\geq 0} \frac{\lambda^k}{k!}(\sup_{t\in [0,T]} |Y_t|)^k \1_{\sup_{t\in [0,T]} |Y_t|\geq 1}\right)\\
&\leq e^\lambda +\E \left( \sum_{k\geq 0} \frac{\lambda^k}{k!}(\sup_{t\in [0,T]} |Y_t|)^{2k}\right).\\
\end{align*}
Using Proposition \ref{estimeesclassiques}  (see the Appendix), this latter term is dominated by
$$e^\lambda +\E \left( \sum_{k\geq 0} \frac{\lambda^k}{k!}(C(1+T^k)\left\{\E(|Y_0|^{2k})+\left(1+\E(|Y_0|^{2k})\right)T^ke^{CT}\right\}\right),$$
which is finite since $Y_0^2$ has a finite Laplace transform. The same reasoning can be applied to $Y_t^n$ which satisfies (\ref{eulerpourY}).
\hfill \framebox[0.6em]{}

\bigskip

\noindent In the proof of Theorem \ref{rateofconvergence} we use the following lemma which enables us to approximate a given function by a sequence of step functions, instead of an arbitrary sequence, without changing the rate of convergence. This lemma is proved in the Appendix.

\begin{lem}
\label{lem:1}
Let $f$ and $(f_n)_n$ denote functions defined on $[0, T]$. For all integer $n\geq 1$, define $\bar{f}_n(t):=f_n(kT/n),\, \, \forall t\in[kT/n, (k+1)T/n]$. Let $\ell>0$, and suppose that $f$ is H\" older continuous with exponent $\beta$ for every $\beta<\ell$. Then,
\begin{equation*}
\forall\, \,  \alpha<\ell,\, \,  \lim_{n\rightarrow +\infty}  \, \, \, n^\alpha \sup_{t\in[0,T]} |f_n(t)-f(t)|=0
\end{equation*}
is equivalent to
\begin{equation*}
\forall\, \,  \alpha<\ell,\, \, \lim_{n\rightarrow +\infty}  \, \, \, n^\alpha \sup_{t\in[0,T]} |\bar{f}_n(t)-f(t)|=0.
\end{equation*}
\end{lem}

\noindent\textbf{Proof of Theorem \ref{rateofconvergence}.}
We prove the statement in detail for ${\tilde Z}_n$ (the proof for $Z_n$ is similar). This is equivalent to prove that 
\begin{equation*}
\forall \, \, 0\leq \alpha < \frac{1}{4}, \,  \, \, \,  \lim_{n\rightarrow +\infty}  \,  n^\alpha \, sup_{t\in [0,T]} {|\bar{X}^n(|\bar{Y}_t^n|)-X(|Y_t|)|}=0,\, \,  {\rm a.s.},
\end{equation*}
thanks to Lemma \ref{lem:1} applied to $Z_t=X(|Y_t|)$ which is locally H\" older continuous with exponent $\beta$ for all $0<\beta<1/4$.

\noindent Let $\alpha\in]0,1/4[$. Let us take $\epsilon>0$ and $M'>0$ such that $\P(\sup_{t\in [0,T]} |Y_t|\geq M')\leq \epsilon$ and set $A:=\{\sup_{t\in [0,T]} |Y_t|<M'\}$. Since $\bar{Y}^n$ converges a.s. uniformly to $Y$ on $[0,T]$, there exists $N$ such that a.s., for all $n\geq N$, $\sup_{t\in [0,T]} |\bar{Y}_t^n-Y_t|<1$. In particular the inequality $\sup_{t\in [0,T]} |\bar{Y}_t^n|<1+M'$ holds a.s. on $A$. Since $X$ is H\" older continuous on $[0,M'+1]$ of exponent $\beta$ for all $\beta\in]2\alpha,\frac{1}{2}[$, there exists a random value $C$ such that for all $n\geq N$, a.s. on $A$,
\begin{equation*}
|X_{|Y_t|}-\bar{X}^n_{|\bar{Y}_t^n|}|\leq C|Y_t-\bar{Y}_t^n|^\beta+ |X_{|\bar{Y}_t^n|}-\bar{X}^n_{|\bar{Y}_t^n|}|.
\end{equation*}
Taking the supremum and multiplying by $n^\alpha$, we obtain that for all $n\geq N$, a.s. on $A$, 
\begin{equation*}
n^\alpha \sup_{[0,T]} |X_{|Y_t|}-\bar{X}^n_{|\bar{Y}_t^n|}|\leq C(n^{\alpha/\beta} \sup_{[0,T]}|Y_t-\bar{Y}_t^n|)^\beta
+n^\alpha \sup_{[0,M'+1]}|X_t-\bar{X}_t^n|).
\end{equation*}
It is proved in \cite{faure} that
\begin{equation*}\label{cv1}
\forall \, \, 0\leq \alpha < \frac{1}{2}, \,  \, \, \,  \lim_{n\rightarrow +\infty}  \,  n^\alpha sup_{t\in [0,T]} {|\tilde{X}_t^n-X_t|}=0, \, \,  {\rm a.s.}
\end{equation*}
Lemma \ref{lem:1} applied to $X$ implies that the same convergence holds for $(\bar{X}^n)$ instead of $\tilde{X}^n$ namely
\begin{equation*}\label{cv2}
\forall \, \, 0\leq \alpha < \frac{1}{2}, \,  \, \, \,  \lim_{n\rightarrow +\infty}  \,  n^\alpha sup_{t\in [0,T]} {|\bar{X}_t^n-X_t|}=0, \, \,  {\rm a.s.}
\end{equation*}
The same argument applies to $Y$. This leads to 
\begin{equation*}
\lim_{n\rightarrow +\infty}  \, n^\alpha \sup_{t\in [0,T]} |X_{|Y_t|}-\bar{X}^n_{|\bar{Y}_t^n|}|=0.
\end{equation*}
Therefore
\begin{equation*}
\P\left(\lim_{n\rightarrow +\infty}  \, n^\alpha \sup_{t\in [0,T]} |X_{|Y_t|}-\bar{X}^n_{|\bar{Y}_t^n|}|=0\right)\geq \P(A)>1-\epsilon.
\end{equation*}
We conclude the proof by letting $\epsilon$ tend to $0$.
\hfill \framebox[0.6em]{}

\bigskip

\section{PDEs for iterated processes with general time and position processes}\label{pdes}

\noindent In this section, we connect iterated processes and high order PDEs in the spirit of \cite{allouba}, \cite{erkan}. However we consider a general framework where the position or time process is not necessarily Markovian. It seems difficult to prove a general statement when the density $p_Y$ of the time process satisfies a PDE with space dependent coefficients. In this case we treat an example. In particular we address the iteration by a Brownian motion with drift and by an Ornstein-Uhlenbeck process as time process. When the time process is a Brownian motion or an $\alpha$-stable process, the equation we obtain reduces to the equation obtained respectively by \cite{allouba} and \cite{erkan}. However our result holds under more general assumptions on the initial value. Moreover in this particular case we prove that the solution of the PDE is unique.

\bigskip

\subsection{Iteration by a general time process}

\noindent In this section the time process $(Y_t)_{t\geq 0}$ is real valued, starts from $0$ at time $0$ and satisfies the following assumptions. 

\medskip

\noindent (A1) $Y_t$ admits a density denoted below by $p_Y(t,0,u)$ (or simply $p_Y$) for all $t>0$,

\smallskip

\noindent (A2) for all $t>0$, $p_Y(t,0,\cdot)\in \mathcal{C}^{r}(\R)$ and $\partial^i_t \partial^j_u \, p_Y(t,0,\cdot)\in L^1(\R)$ for all  $0\leq i\leq q,\ 0\leq j \leq r$, for some integers $q,r$, 

\smallskip

\noindent (A3) the even part $\mathcal{E}(p_Y)$ of $p_Y$ satisfies
\begin{equation}
\label{edpdensiteY}
\sum_{k=1}^{p} \alpha_k \frac{\partial^k}{\partial t^k}\mathcal{E}(p_Y)(t,0,u)=\sum_{i=0}^{q} \sum_{j=1}^{r} \beta_{i,j}\frac{\partial^{i+j}}{\partial t^i \partial u^j}\mathcal{E}(p_Y)(t,0,u),
\end{equation}
 for some integer $p$ and real valued functions of the time variable $\alpha_k$, $\beta_{i,j}$ (integers $q,r$ are those appearing in (A2)).

\bigskip

\noindent Examples satisfying (A1)-(A3) are provided by 
\begin{itemize}
\item Brownian motion corresponding to $p=1,\, q=0,\, r=2,\,$ $\alpha_1=1,\, \beta_{0,1}=0,\, \beta_{0,2}=1/2$, 
\item the Cauchy process $p=2,\, q=0,\, r=2,\,$$\alpha_1=0,\, \alpha_2=1,\, $ $\beta_{0,1}=0, \, \beta_{0,2}=-1$, 
\item the telegraph process whose density satisfies
\begin{equation*}
\left( \frac{\partial^2}{\partial t^2}+2\lambda \frac{\partial}{\partial t}\right)p_Y(t,0,u)=v^2\frac{\partial^2}{\partial u^2}p_Y(t,0,u), \quad \lambda, v>0.
\end{equation*}
\item the sum of a telegraph process and an independent Brownian motion. In this case $p_Y$ satisfies (see \cite{Blanchard1993225})
\begin{multline*}
\left( \frac{\partial^2}{\partial t^2}+2\lambda \frac{\partial}{\partial t}\right)p_Y(t,0,u)\\
=\left((v^2+\lambda)\frac{\partial^2}{\partial u^2}+\frac{\partial^3}{\partial t \partial u^2}-\frac{1}{4}\frac{\partial^4}{\partial u^4}\right)p_Y(t,0,u), \quad \lambda, v>0.
\end{multline*}
\end{itemize}
Actually in the first three cases the density $p_Y(t,0,\cdot)$ itself is even for all $t$. 

\medskip

\noindent If $p_Y$ satisfies (\ref{edpdensiteY}) and if moreover the orders of the partial derivatives w.r.t. the space variable $u$ on the right hand side of (\ref{edpdensiteY}) are even, then $\mathcal{E}(p_Y)$ satisfies (\ref{edpdensiteY}) as well.

\medskip

\noindent Another example is obtained when $Y$ is a Brownian motion with variance $\sigma^2$ and constant drift $\mu$. From the Fokker-Planck equation satisfied by $p_Y,$ 
\begin{equation*}
\frac{\partial}{\partial t}p_Y(t,0,u)=\frac{1}{2}\sigma^2\frac{\partial^2}{\partial u^2}p_Y(t,0,u)-\mu\frac{\partial}{\partial u}p_Y(t,0,u),\quad \sigma,\mu\in \R,
\end{equation*}
we deduce that its even and odd parts satisfy the system 
\begin{align*}
\frac{\partial}{\partial t}\mathcal{E}(p_Y)&=\frac{1}{2}\sigma^2\frac{\partial^2}{\partial u^2}\mathcal{E}(p_Y)-\mu\frac{\partial}{\partial u}\mathcal{O}(p_Y),\\
\frac{\partial}{\partial t}\mathcal{O}(p_Y)&=\frac{1}{2}\sigma^2\frac{\partial^2}{\partial u^2}\mathcal{O}(p_Y)-\mu\frac{\partial}{\partial u}\mathcal{E}(p_Y).
\end{align*}
which implies the following particular case of (\ref{edpdensiteY}) 
\begin{equation*}
\frac{\partial^2}{\partial t^2}\mathcal{E}(p_Y)-\sigma^2\frac{\partial^3}{\partial t\partial u^2}\mathcal{E}(p_Y)-\mu^2\frac{\partial^2}{\partial u^2}\mathcal{E}(p_Y)+\frac{1}{4}\sigma^4\frac{\partial^4}{\partial u^4}\mathcal{E}(p_Y)=0.
\end{equation*}

\bigskip

\bigskip

\noindent  In this subsection the position process is a Markov process $(X_t^x)_{t\geq 0}$ starting from a given $x\in \R$ with semigroup $P^t$ and infinitesimal generator $\mathcal{L}$ with domain $D(\mathcal{L})$. We denote by $B(\R^d,\R^d)$ the set of bounded functions from $\R^d$ to $\R^d$ and define $B_0:=\{f \in B(\R^d,\R^d), \|P^tf-f\|_\infty \rightarrow 0 \text{\ as\ } t\downarrow 0\}$. $B_0$ is a closed subset of $B(\R^d,\R^d)$ containing $D(\mathcal{L})$. In the following, we will use iterations of $\mathcal{L}$:
\begin{equation*}
D(\mathcal{L}^0)=B_0 \text{\quad and\quad}
D(\mathcal{L}^n)=\{f\in D(\mathcal{L}^{n-1});\,  \mathcal{L}^{n-1}f\in D(\mathcal{L})\},\, \, \forall\, n\geq 1.
\end{equation*}

\noindent We state below the results of this section (Theorem \ref{FKitere} and its corollaries). Their proofs are provided in section \ref{proofs}. 

\begin{thm}
\label{FKitere}
Let $Y$ be a real valued process independent of $X^x$ whose density $p_Y$ satisfies assumptions {\rm (A1)}-{\rm (A3)}. Let $f\in D(\mathcal{L}^{r-1})$ where $r$ is the highest order of the partial derivatives in $u$ appearing in (\ref{edpdensiteY}).  Define $v(t,x):=\E \left[f(X^x(|Y_t|))\right]$. Then $v\in D(\mathcal{L}^{r})$ and satisfies
\begin{equation}
\label{edpFKitere}
\sum_{k=1}^{p} \alpha_k \frac{\partial^k}{\partial t^k}v(t,x)=\sum_{i=0}^{q}\sum_{j=1}^{r}\beta_{i,j}(-1)^j\frac{\partial^i}{\partial t^i}\mathcal{L}^jv(t,x)+{\cal B}_f(t,x),
\end{equation}
on $]0,+\infty[\times \R^d$ where $\mathcal{L}$ acts on $x$ and  ${\cal B}_f(t,x)$ is the boundary term
\begin{equation}
\sum_{i=0}^{q}\sum_{j=1}^{r}\, \, \beta_{i,j}\sum_{0\leq 2k\leq j-1} (-1)^{j-2k}\mathcal{L}^{j-1-2k}f(x)\left. \frac{\partial^{i+2k}}{\partial t^i \partial u^{2k}}p_Y(t,0,u)\right|_{u=0}.
\end{equation}
\end{thm}

\bigskip

\noindent Restricting to $Y$ a Brownian motion, we recover as a corollary the following result of \cite{allouba}. Moreover we show that in this case (\ref{edpFKitere}) admits a unique solution under weaker assumptions on the initial condition.
\begin{corol}\label{corol:2}
Let $f\in D(\mathcal{L})$ and $B$ a Brownian motion independent of $X^x$. Then, $v(t,x):=\E \left[f(X^x(|B_t|))\right]$ is an element of $D(\mathcal{L}^2)$ and is the unique solution in $D(\mathcal{L}^2)$ of
\begin{equation}
\begin{cases} 
\label{FKibm}
\frac{\partial}{\partial t}v(t,x)=\frac{1}{\sqrt{2\pi t}}\mathcal{L}f(x)+\frac{1}{2}\mathcal{L}^2v(t,x), \quad t>0,\ x\in\R^d,\\
v(0,x)=f(x), \quad x\in\R^d.
\end{cases}
\end{equation}

\end{corol}
\begin{corol}
Let $f\in D(\mathcal{L}^3)$ and $\, Y$ be a Brownian motion with drift $\mu$ and diffusion coefficient $\sigma$, independent of $X^x$. Then, there exist two functions of time, $\alpha$ and $\beta$, such that $v(t,x):=\E \left[f(X^x(|Y_t|))\right]$ satisfies
\begin{multline*}
\frac{\partial^2}{\partial t^2}v(t,x)=\sigma^2\frac{\partial}{\partial t}\mathcal{L}^2v(t,x)+\mu^2\frac{\partial}{\partial t}\mathcal{L}^2v(t,x)-\frac{1}{4}\sigma^4\frac{\partial}{\partial t}\mathcal{L}^4v(t,x)\\
+\alpha(t)\mathcal{L}f(x)+\beta(t)\mathcal{L}^3f(x), \quad \forall \, t>0, \forall \, u\in \R,
\end{multline*}
with $v(t,\cdot)\in D(\mathcal{L}^4)$, $\forall\,  t>0$. By definition $v(0,x)=f(x)$. 
\end{corol}

\noindent The time processes considered in the above statements are associated to PDEs whose coefficients may depend on the time variable but not on the spatial variable. It seems difficult to extend Theorem \ref{FKitere} to a large class of time processes. For instance, if $Y$ is an Ornstein-Uhlenbeck process issued from $0$ satisfying $dY_t=-\frac{1}{2}Y_tdt+dW_t$, its density satisfies the Fokker-Planck equation 
\begin{equation}
\label{FPOU}
\frac{\partial}{\partial t}p_Y(t,0,u)=\frac{1}{2}\frac{\partial^2}{\partial u^2}p_Y(t,0,u)+\frac{1}{2}\frac{\partial}{\partial u}\left[up_Y(t,0,u)\right], \,  \forall\,  t>0, \forall\,  u\in \R.
\end{equation}

\noindent Since some coefficients of this equation are functions of the spatial variable $u$, we cannot apply Theorem \ref{FKitere} directly. However in this particular case a rewriting of the equation leads to the following result. The difficulty is that in general such a rewriting is not possible.

\begin{prop}\label{OU}
Let $X^x$ be a $\R^d$-valued Markov process with infinitesimal generator $\mathcal{L}$. Let $Y$ an Ornstein-Uhlenbeck process independent of $X$, satisfying $dY_t=-\frac{1}{2}Y_tdt+dW_t$, $Y_0=0$. Then, for $f\in D(\mathcal{L})$, the function $v(t,x):=\E \left[f(X^x(|Y_t|))\right]$ satisfies $v(t,\cdot)\in D(\mathcal{L}^2)$ for all $t>0$ and solves
\begin{equation}
\label{eq:46}
\frac{\partial}{\partial t}v(t,x)=\frac{e^{-t}}{2\sqrt{2\pi(1-e^{-t})}}\mathcal{L}f+\frac{1}{2}e^{-t}\mathcal{L}^2v(t,x), \quad t>0,\ x\in\R^d,
\end{equation}
with initial condition $v(0,x)=f(x)$.
\end{prop}

\begin{corol}
Let $Y$ be a telegraph process with parameters $\lambda>0$ and $v>0$, independent of $X$ and $f\in D(\mathcal{L})$. Then, $u(t,x)=\E \left[f(X^x(|Y_t|))\right]$ is solution of
\begin{equation*}
\begin{cases} 
\left(\frac{\partial^2}{\partial t^2}+\lambda\frac{\partial}{\partial t}\right)u(t,x)=v^2p_Y(t,0,0)\mathcal{L}f(x)+v^2\mathcal{L}^2u(t,x), \quad t>0,\ x\in\R^d.\\
u(0,x)=f(x), \quad x\in\R^d
\end{cases}
\end{equation*}
\end{corol}

\subsection{Proofs}\label{proofs}
\noindent We now come to the proofs of Theorem \ref{FKitere}, Corollary \ref{corol:2} and Proposition \ref{OU}. We start with the proof of Theorem \ref{FKitere} which relies on the following lemma proved in the Appendix (section \ref{appendix}).

\begin{lem}
\label{IPP}
Let $g\in L^1(\R)$ such that $g$ is differentiable with $g'\in L^1(\R)$. Then for all $f\in B_0$, the function 
$F(g):x\mapsto \int_0^\infty g(s)P^sf(x)ds$ belongs to $D(\mathcal{L})$ and
$$\mathcal{L}F(g)=-F(g')-g(0)f.$$
More generally, if $f\in D(\mathcal{L}^{n-1})$ and for all $k\in\{0,\ldots, n\}$, $g^{(k)}\in L^1(\R)$,  then $F(g)\in D(\mathcal{L}^n)$ and
\begin{equation}
\mathcal{L}^nF(g)=(-1)^nF(g^{(n)})-\sum_{l=0}^{n-1}(-1)^{n-1-l}g^{(n-1-l)}(0)\mathcal{L}^lf. \label{formuleIPP}
\end{equation}
\end{lem}

\medskip

\noindent\textbf{Proof of Theorem \ref{FKitere}.}
The independence of $X$ and $Y$ implies that
\begin{equation*}
v(t,x)=\int_\R p_Y(t,0,u)P^{|u|}f(x)du,
\end{equation*}
and the fact that $\mathcal{E}(p_Y)(t,0,u)=\frac{1}{2}(p_Y(t,0,u)+p_Y(t,0,-u))$ yield
\begin{equation}
v(t,x)=2\int_0^\infty \mathcal{E}(p_Y)(t,0,u)P^uf(x)du. 
\end{equation}
For fixed $t>0$ define $g(\cdot):=\mathcal{E}(p_Y)(t,0,\cdot)$. Then $v(t,\cdot)$ coincides with $2F(g)$. From {\bf (A3)} we obtain
\begin{eqnarray*}
\sum_{k=1}^{p} \alpha_k(t)\frac{\partial^k}{\partial t^k}v(t,x)&=&2\int_0^\infty \sum_{k=1}^{p} \alpha_k(t)\frac{\partial^k}{\partial t^k}\mathcal{E}(p_Y)(t,0,u)P^u f(x)du\nonumber \\
&=&2\sum_{i=0}^{q}\sum_{j=1}^{r} \beta_{i,j}(t)\frac{\partial^i}{\partial t^i}\int_0^\infty \frac{\partial^j}{\partial u^j}\mathcal{E}(p_Y)(t,0,u)P^u f(x)du\\
&=& 2\sum_{i=0}^{q}\sum_{j=1}^{r} \beta_{i,j}(t)\frac{\partial^i}{\partial t^i} F(g^{(j)})(x),\label{intermediaire}
\end{eqnarray*}
where we make explicit the fact that $\alpha_k$ and $\beta_{i,j}$ depend only on time. Applying (\ref{formuleIPP}) to each $0\leq j\leq r$, we obtain
\begin{equation*}
F(g^{(j)})=(-1)^j{\cal L}^j F(g)+\sum_{\ell=0}^{j-1}(-1)^{\ell+1}\left. \frac{\partial^{j-1-\ell}}{\partial u^{j-1-\ell}}\mathcal{E}(p_Y)(t,0,u)\right|_{u=0}\mathcal{L}^\ell f.
\end{equation*}
We can now conclude using the fact the even function $g\equiv \mathcal{E}(p_Y)(t,0,\cdot)$ satisfies 
\begin{equation*}
\left. \frac{\partial^\ell}{\partial u^\ell}\mathcal{E}(p_Y)(t,0,u)\right|_{u=0}=
\left\lbrace 
\begin{array}{lcl}
\left. \frac{\partial^\ell}{\partial u^\ell}p_Y(t,0,u)\right|_{u=0} &\text{if $\ell$ even}\\ 
0 & \text{if $\ell$ odd}\\
\end{array}\right.
\end{equation*}
for all $0\leq \ell\leq j-1$. 
\hfill \framebox[0.6em]\\

\bigskip

\noindent {\bf Proof of Corollary \ref{corol:2}}
We deduce from Theorem \ref{FKitere} that $v\in D(\mathcal{L}^2)$ and it satisfies (\ref{FKibm}). Let $\tilde{v}$ be another solution of  (\ref{FKibm}) belonging to $D(\mathcal{L}^2)$. Then $u=v-\tilde{v}\in D(\mathcal{L}^2)$ and $u$ is solution of
\begin{equation}
\label{eq:44}
\frac{\partial}{\partial t}u(t,x)=\frac{1}{2}\mathcal{L}^2u(t,x), \quad t>0,\ x\in\R,
\end{equation}
with initial condition $u(0,x)=0, \forall x\in \R^d$. Since $\mathcal{L}^2$ generates a strongly continuous semigroup on $B_0$ (see \cite{nag86}), the solution of equation (\ref{eq:44}) is unique and vanishes identically, which implies that $v\equiv \tilde{v}$.\hfill \framebox[0.6em]

\bigskip

\noindent{\bf Proof of Proposition \ref{OU}}
Since $p_Y(t,0,u)=\frac{1}{\sqrt{2\pi(1-e^{-t})}}e^{-\frac{u^2}{2(1-e^{-t})}}$, it satisfies the two following identities,
\begin{equation*}
u\,  p_Y(t,0,u)=-(1-e^{-t})\frac{\partial}{\partial u}p_Y(t,0,u),
\end{equation*}
\begin{equation}\label{pdeOU}
\frac{\partial}{\partial t}p_Y(t,0,u)=\frac{1}{2}e^{-t}\frac{\partial^2}{\partial u^2}p_Y(t,0,u), \quad \forall\,  t>0,\, \forall\,  u\in \R.
\end{equation}
The argument in the proof of Theorem \ref{FKitere} applies to (\ref{pdeOU}) instead of the original equation (\ref{FPOU}). Thus (\ref{eq:46}) is satisfied. (\ref{eq:46}) is indeed a particular case of (\ref{edpFKitere}) where $p=1, q=0, r=2$, $\beta_{0,1}=0,\, \beta_{0,2}=\frac{1}{2}e^{-t}, \alpha_1=1$. The boundary term $\mathcal{B}_f(t,x)$ is equal to $\frac{1}{2}e^{-t}p_Y(t,0,0)=\frac{e^{-t}}{2\sqrt{2\pi(1-e^{-t})}}.$ \hfill \framebox[0.6em]{}

\begin{rmq}
Let $r:=\inf \{k\in \N,\, \partial^k_u \, p_Y(t,0,u)|_{u=0}\not \equiv 0\}$. Then the assumption $f\in D(\mathcal{L}^{n_j-1})$ can be weakened and replaced by $f\in D(\mathcal{L}^{n_j-1-r})$ by a slight modification of equality (\ref{intermediaire}). For instance, if $0\not \in supp\, p_Y(t,0,\cdot) $, for all $t>0$, , then we only need $f\in B_0$ and equation (\ref{edpFKitere}) becomes
\begin{equation*}
\sum_{k=1}^p \alpha_k \frac{\partial^k}{\partial t^k}v(t,x)=\sum_{i=0}^{q}\sum_{j=1}^{r} \beta_{i,j}(-1)^j.\frac{\partial^i}{\partial t^i}\mathcal{L}^jv(t,x).
\end{equation*}
Obviously if $Y$ is a non-negative process and $p_Y$ satisfies an equation of the type (\ref{edpdensiteY}), then removing the absolute value, $v(t,x):=\E \left[f(X^x(Y_t))\right]$ is solution of (\ref{edpFKitere}).
\end{rmq}

\subsection{General position processes}\label{gyongy}

In this subsection, the position process $X^x$ is not necessarily markovian. Neither is the time process for which we assume {\rm (A1)}-{\rm (A3)}. The position process $X^x$ is an It\^o process, it can be written as
\begin{equation}
X_t^x=x+\int_0^t \delta(s,\omega)dW_s+\int_0^t \beta(s,\omega)ds,
\end{equation}
where $(W(t),\mathcal{F}_t)$ is a Wiener process, $\delta$ and $\beta$ are bounded $\mathcal{F}_t$-nonanticipative processes such that $\delta\delta^T$ is uniformly positive definite.

\begin{thm}
Let us keep the assumptions of Theorem \ref{FKitere} on $Y$. Let $\mathcal{L}$ be a differential operator defined by
$$\mathcal{L}:=\frac{1}{2}\sigma(x)^2\frac{\partial^2}{x^2}+b(x)\frac{\partial}{\partial x}$$
with
\begin{eqnarray}
\sigma(t,z)&:=&\E(\delta \delta^T(t)\, |\, \ X_t^x=z) \label{eq:95}\\
b(t,z)&:=&\E(\beta(t)\, |\, \ X_t^x=z). \label{eq:96}
\end{eqnarray}
Let $f\in D(\mathcal{L}^{r-1})$ where $r$ is the highest order of the partial derivatives in $u$ appearing in (\ref{edpdensiteY}).  Define $v(t,x):=\E \left[f(X^x(|Y_t|))\right]$. Then $v\in D(\mathcal{L}^{r})$ and satisfies
\begin{equation}
\label{edpFKitere2}
\sum_{k=1}^{p} \alpha_k \frac{\partial^k}{\partial t^k}v(t,x)=\sum_{i=0}^{q}\sum_{j=1}^{r}\beta_{i,j}(-1)^j\frac{\partial^i}{\partial t^i}\mathcal{L}^jv(t,x)+{\cal B}_f(t,x),
\end{equation}
on $]0,+\infty[\times \R^d$ where $\mathcal{L}$ acts on $x$ and  ${\cal B}_f(t,x)$ is the boundary term
\begin{equation}
\sum_{i=0}^{q}\sum_{j=1}^{r}\, \, \beta_{i,j}\sum_{0\leq 2k\leq j-1} (-1)^{j-2k}\mathcal{L}^{j-1-2k}f(x)\left. \frac{\partial^{i+2k}}{\partial t^i \partial u^{2k}}p_Y(t,0,u)\right|_{u=0}.
\end{equation}
\end{thm}

\begin{dem}
Under assumptions on $\delta$ and $\beta$, there exists from \cite{mimicking} a process $({\cal X}_t)$ weak solution of
\begin{equation}
{\cal X}_t^x=x+\int_0^t \sigma(s,{\cal X}_s)dW_s+\int_0^t b(s,{\cal X}_s)ds,
\end{equation} 
such that $\forall\, t>0$, ${\cal X}_t^x$ and $X_t^x$ are identically distributed and coefficients $b$ and $\sigma$ are given by equations (\ref{eq:95}) and (\ref{eq:96}).\\
Then for all $t$, $\E(f(X_t^x))=\E(f({\cal X}^x_t))$. For $(Y_t)$ independent of $X$ and of ${\cal X}$, we have 
\begin{equation}
v(t,x):=\E \left[f(X^x(|Y_t|))\right]=\E \left[f({\cal X}^x(|Y_t|))\right]
\end{equation}
therefore
\begin{equation}\label{newposition}
v(t,x):=\E \left[f(X^x(|Y_t|))\right]=\E \left[f({\cal X}^x(|Y_t|))\right]=2\int_0^\infty \mathcal{E}(p_Y)(t,0,u)P^uf(x)du, 
\end{equation}
where $P^u$ in (\ref{newposition}) denotes the semigroup of ${\cal X}$ at time $u$. Hence under assumptions {\bf (A1)}-{\bf (A3)}, the result of Theorem \ref{FKitere} applies to ${\cal X}$, $Y$ and the function $v$ which satisfies PDE (\ref{edpFKitere2}) with ${\cal L}$ the infinitesimal generator of ${\cal X}$.\\

\end{dem}

\section{Transformations of high order PDEs}\label{transformations}

\noindent In this section we consider transformations between two high order PDEs. We start with the use of intertwining diffusions to build such a transformation.

\subsection{Intertwining diffusions}\label{sec:intertwining}

\noindent Consider two diffusions $(X_t)_{t\geq 0}$ and $(U_t)_{t\geq 0}$ which take values in $\R^d$ and $ \R^n$ respectively, with respective semigroups $P_t^X$ and  $P_t^U$. The diffusions $X$ and $U$ are intertwining if there exists a density function $\Lambda$ such that the operator $L$ defined by
\begin{equation}
L\phi(z):=\int_{\R^d} \Lambda(z,x)\phi(x)dx, \quad \forall\, z\in\R^n,
\end{equation}
satisfies $P_t^U\, L=L\, P_t^X$ for all $t>0$. 

\bigskip
\noindent Assume moreover that there exists a $\R^d\times \R^n$-valued diffusion process $(Z_1,Z_2)$ satisfying:

\noindent (i) $X=Z_1$ and $U=Z_2$ in law,

\noindent (ii) $\E(f(Z_1(0))|Z_2(0)=z_2)=Lf(z_2)$,

\noindent (iii) for any $t>0$ the random variables $Z_2(0)$ and $Z_1(t)$ are conditionally independent given $Z_1(0)$,

\noindent (iv) for any $t>0$ the random variables $Z_2(0)$ and $Z_1(t)$ are conditionally independent given $Z_2(t)$.

\noindent In this case it is proved in \cite{pal} that $\Lambda$ is a distributional solution of the hyperbolic PDE
\begin{equation}\label{edppourLambda}
({\cal L}_X)^* \; \Lambda={\cal L}_U\, \Lambda,
\end{equation}
where $\mathcal{L}_X$ and $\mathcal{L}_U$ denote the infinitesimal generators of $X$ and $U$.

\bigskip

\noindent In the following theorem we show that intertwining preserves the structure of (\ref{FKitere}), provided that we iterate $X$ and $U$ by the same time process. Let us recall that if $X^x$, $Y$ and $f$ satisfy the assumptions of Theorem \ref{FKitere}, then $v(t,x):=\E \left[f(X^x(|Y_t|))\right]$ is a solution of (\ref{edpFKitere}) with initial condition $v(0,\cdot)\equiv f(\cdot)$. We rewrite (\ref{edpFKitere}) as 
\begin{equation}
\sum_{k=1}^{p} \alpha_k \frac{\partial^k}{\partial t^k}v(t,x)=Q(\frac{\partial}{\partial t},\mathcal{L}_X)v(t,x)+\mathcal{B}_{f,\mathcal{L}^X}(t,x),
\end{equation}
where the polynomial $Q(\frac{\partial}{\partial t},\mathcal{L}_X)$ coincides with $\sum_{i=0}^{q}\sum_{j=1}^{r}\beta_{i,j}(-1)^j\frac{\partial^i}{\partial t^i}\mathcal{L}_X^j$.

\begin{thm}\label{intertwining}
Let $X^x$ be a diffusion. Assume that $X^x, Y$ and $f$ satisfy the assumptions of Theorem \ref{FKitere}. Let $U^x$ be another diffusion independent of $Y$ such that $X$ and $U$ are intertwining. Define $g:=Lf$ and $h(t,x):=\E(g(U^x_{|Y_t|}))$. Then $h$ belongs to $\mathcal{D}((\mathcal{L}_U)^{r-1})$ ($r$ is the highest order of the partial derivatives in $u$ appearing in (\ref{edpdensiteY})) and satisfies the PDE
\begin{equation}
\label{eq:83}
\sum_{k=1}^{p} \alpha_k \frac{\partial^k}{\partial t^k}h(t,x)=Q(\frac{\partial}{\partial t},\mathcal{L}_U)h(t,x)+\mathcal{B}_{g,\mathcal{L}_U}(t,x),
\end{equation}
on $]0,+\infty[\times \R^n$ with initial condition $h(0,\cdot)\equiv g(\cdot)$ by definition.
\end{thm}
\begin{rmq}
An interesting point is that this theorem maps a solution defined on $\R_+\times \R^d$ into a solution defined on $\R_+\times \R^n$. A version without the boundary term $\mathcal{B}_{f}$ can be derived from Theorem \ref{edpsansf} by a slight modification of the proof. 
\end{rmq}
\bigskip

\noindent {\bf Proof of Theorem \ref{intertwining}.} 
From the relation $P_t^UL=LP_t^X$, we have for all $w\in\mathcal{D}(\mathcal{L}_X)$,
$$\frac{P_t^ULw-Lw}{t}=\frac{LP_t^Xw-Lw}{t}$$
and using the boundedness of $L$, the right-hand side converges. Therefore so does the left-hand side which implies that $Lw\in\mathcal{D}(\mathcal{L}_U)$. Hence, $f\in\mathcal{D}((\mathcal{L}_X)^{r-1})$ implies $Lf\in\mathcal{D}((\mathcal{L}_U)^{r-1})$. Moreover by definition 
\begin{equation*}
h(t,x)=\int_\R p_Y(t,0,\tau)\E(Lf(U_{|\tau|}^x))d\tau=\int_\R p_Y(t,0,\tau) \, P_{|\tau|}^U  Lf(x) \, d\tau.
\end{equation*}
Using first that $P^U L=LP^X$, and then the kernel $\Lambda$, we obtain 
\begin{eqnarray*}
h(t,x)&=&\int_\R p_Y(t,0,\tau) \, LP_{|\tau|}^X  f(x) \, d\tau\\
&=&\int_\R p_Y(t,0,\tau) \, \int_{\R^d} \Lambda (x, \rho) P_{|\tau|}^X  f(\rho) d\rho \, d\tau.
\end{eqnarray*}
We can apply Fubini Theorem which yields
\begin{equation*}
h(t,x)=\int_{\R^d} \Lambda (x, \rho) \int_\R p_Y(t,0,\tau) \, P_{|\tau|}^X  f(\rho)  \, d\tau d\rho,
\end{equation*}
and by definition of $v$, we conclude that
\begin{equation}\label{relationhv}
h(t,\cdot)=\int_{\R^d} \Lambda (\cdot, \rho) v(t,\rho) d\rho.
\end{equation}
We apply $\sum_{k=1}^{p} \alpha_k \frac{\partial^k}{\partial t^k}$ to (\ref{relationhv}), exchange the order of $\frac{\partial^k}{\partial t^k}$ and the integral, and use (\ref{eq:83}) to obtain
\begin{equation}
\label{eq:50}
\sum_{k=1}^{p} \alpha_k \frac{\partial^k}{\partial t^k} h(t,\cdot)=\int_{\R^d} \Lambda (\cdot, \rho) \left\{Q(\frac{\partial}{\partial t},\mathcal{L}_X)v(t,\rho)+\mathcal{B}_{f,\mathcal{L}_X}(t,\rho)\right\} \, d\rho.
\end{equation}
In order to simplify the proof, we consider $Q(M,N)=M^iN^j$ and $\mathcal{B}_{f,\mathcal{L}_X}(t,\rho)=(\mathcal{L}_X)^kf(\rho)\alpha(t)$. We now decompose the right hand side of (\ref{eq:50}) into two parts.
\begin{multline*}
\int_{\R^d} \Lambda (\cdot, \rho) Q(\frac{\partial}{\partial t},\mathcal{L}_X)v(t,\rho) \, d\rho=\int_{\R^d} \Lambda (\cdot, \rho) \frac{\partial^i}{\partial t^i}(\mathcal{L}_X)^jv(t,\rho) \, d\rho\\
=\frac{\partial^i}{\partial t^i}\int_{\R^d} \Lambda (\cdot, \rho) (\mathcal{L}_X)^jv(t,\rho) \, d\rho.
\end{multline*}
Taking the adjoint brings the term $({\cal L}_X)^*\Lambda(\cdot,\rho)$ under the integral. By assumption (\ref{edppourLambda}), applied $j$-times, we conclude that 
\begin{equation*}
\frac{\partial^i}{\partial t^i}\int_{\R^d} \Lambda (\cdot, \rho) (\mathcal{L}_X)^jv(t,\rho) \, d\rho=
\frac{\partial^i}{\partial t^i}\int_{\R^d} (\mathcal{L}_U)^j\Lambda (\cdot, \rho) v(t,\rho) \, d\rho
\end{equation*}
Inverting the order between $(\mathcal{L}_U)^j$ and $\int$, we obtain
\begin{equation*}
\frac{\partial^i}{\partial t^i}\int_{\R^d} \Lambda (\cdot, \rho) (\mathcal{L}_X)^jv(t,\rho) \, d\rho=
\frac{\partial^i}{\partial t^i}(\mathcal{L}_U)^j\int_{\R^d} \Lambda (\cdot, \rho) v(t,\rho) \, d\rho
\end{equation*}
So by definition of $h$, we have shown
\begin{equation}
\label{eq:51}
\int_{\R^d} \Lambda (\cdot, \rho) Q(\frac{\partial}{\partial t},\mathcal{L}_X)v(t,\rho) \, d\rho=
Q(\frac{\partial}{\partial t},\mathcal{L}_U)h(t,\cdot)
\end{equation}

\begin{multline*}
\int_{\R^d} \Lambda (\cdot, \rho) \mathcal{B}_{f,\mathcal{L}_X}(t,\rho) \, d\rho=\int_{\R^d} \Lambda (\cdot, \rho) (\mathcal{L}_X)^kf(\rho)\alpha(t) \, d\rho\\
=\int_{\R^d} (\mathcal{L}_U)^k\Lambda (\cdot, \rho) f(\rho)\alpha(t) \, d\rho
=(\mathcal{L}_U)^k\alpha(t)\int_{\R^d} \Lambda (\cdot, \rho) f(\rho) \, d\rho\\
=(\mathcal{L}_U)^k g(\cdot)\alpha(t)=\mathcal{B}_{g,\mathcal{L}_U}(t,\cdot)
\end{multline*}
This combined with (\ref{eq:50}) and (\ref{eq:51}) gives the announced PDE for $h$ thanks to (\ref{relationhv}). In all the identities above ${\cal L}_X$ acts on $\rho$ whereas ${\cal L}_U$ acts on the other variable represented by a $\cdot$ (as in $(\mathcal{L}_U)^k\Lambda (\cdot, \rho))$). 
 \hfill \framebox[0.6em]{}\\
\bigskip

\noindent As an example let $X=X_\alpha$ and $U=X_\alpha+X_\beta$ with $X_\alpha$ and $X_\beta$ two independent squared-Bessel processes of respective dimension $2\alpha>0$ and $2\beta>0$ (cf. \cite{pal}). Their semigroups are intertwining with 
\begin{equation*}
\Lambda(y,x)=\frac{y^{-1}}{B(a,b)}\left(\frac{x}{y}\right)^{\alpha-1}\left(1-\frac{x}{y}\right)^{\beta-1}\1_{0<x<y}.
\end{equation*}
Here $B$ denotes the Beta function. The infinitesimal generators are given by
\begin{equation*}
\mathcal{L}_X=2x\partial_{xx}+2\alpha\partial_y\quad \text{and}\quad \mathcal{L}_U=2y\partial_{yy}+2(\alpha+\beta)\partial_y.
\end{equation*}
Let $f$ be a real function of class $C^2$ with compact support so that $f\in \mathcal{D}(\mathcal{L}_X)$ and let $W$ be a Brownian motion independent of $X$ and $U$. Then from Corollary \ref{corol:2}, $v(t,x):=\E \left[f(X^x(|W_t|))\right]$ is the unique solution in $\mathcal{D}((\mathcal{L}_X)^2)$ of
\begin{equation}
\begin{cases} 
\frac{\partial}{\partial t}v(t,x)=\frac{1}{\sqrt{2\pi t}}\mathcal{L}_Xf(x)+\frac{1}{2}(\mathcal{L}_X)^2v(t,x), \quad t>0,\ x\in\R\\
v(0,x)=f(x), \quad x\in\R
\end{cases}
\end{equation}
and from Theorem \ref{intertwining}, $h(t,x):=\E \left[g(U^x(|W_t|))\right]$ is the unique solution in $\mathcal{D}((\mathcal{L}_U)^2)$ of
\begin{equation}
\begin{cases} 

\frac{\partial}{\partial t}v(t,x)=\frac{1}{\sqrt{2\pi t}}\mathcal{L}_Ug(x)+\frac{1}{2}(\mathcal{L}_U)^2h(t,x), \quad t>0,\ x\in\R\\
h(0,x)=g(x), \quad x\in\R
\end{cases}
\end{equation}
where $g(x):=Lf(x)=\int_0^1f(x\rho)\rho^{\alpha-1}(1-\rho)^{\beta-1}d\rho.$

\subsection{Mapping a high order PDE into another one}

\medskip

\noindent In this section, we consider two initial-boundary value problems of high order (cf.\cite{bragg}). We will construct a mapping  that transforms one of these problems into the other one. We show that this mapping can be expressed using a Feynman-Kac formula. If  $x=(x_1,\ldots,x_n)$ is a point in ${\R}^n$ and  $\phi$ is a smooth function, we write $D_i\phi(x):=\partial \phi(x)/\partial x_i$. We need compositions of the form $D^\alpha:=D_1^{\alpha_1}D_2^{\alpha_2}\cdots D_k^{\alpha_k}$ where $k$ and the $\alpha_i$ are integers. Finally we denote by $P(x,D)$ any polynomial in the partial derivatives as follows
\begin{equation*}
P(x,D):=\sum_{0\leq\alpha \leq m }a_\alpha(x)D^\alpha,
\end{equation*}
where $a_\alpha(x)$ are given functions of $x$. The first initial-boundary value problem is
\begin{equation}\label{eq:48}
\left\lbrace 
\begin{array}{lcl} 
\partial^2 v(x,t)/\partial t^2=P(x,D)v(x,t), &t>0,\\
v(x,0)=0, \quad v_t(x,0)=\phi(x),\\
B(x,D)v(x,t)=g(x,t),&x\in {\cal S},\, t>0,
\end{array}\right.
\end{equation}
\noindent where $B(x,D)$ is a non tangential boundary operator and ${\cal S}:=\{ x; S(x)=0\}$ for some function $S$, denotes a cylindrical surface. 

\begin{thm}\label{transform}
Suppose that (\ref{eq:48})  admits a solution $v(x,t)$. Let $B$ be a Brownian motion and suppose that the function $u(t,x):=\frac{\partial^{1/2}}{\partial t^{1/2}}\E (v(x,|B_{2t}|))$ is well defined. Then $u(t,x)$ solves the following problem
\begin{equation}\label{eq:49}
\left\lbrace 
\begin{array}{lcl} 
\partial u(x,t)/\partial t=P(x,D)u(x,t), &t>0,\\
u(x,0)=\phi(x),\\
B(x,D)u(x,t)=f(x,t),&x\in S,\, t>0,
\end{array}\right.
\end{equation}
where $f(x,t):=\frac{1}{2\sqrt{\pi}t^{3/2}}\int_0^\infty\xi e^{-\xi^2/4t}g(x,\xi)d\xi$.
\end{thm}
In Theorem \ref{transform} the derivative in $u(t,x)$ is a particular case of the Caputo fractional derivative. Take a function $f$ and a positive real $\gamma$. If $m-1<\gamma<m$ for some integer $m$, the Caputo derivative of $f$ at order $\gamma$ is  

\begin{equation*}
\frac{\partial^\gamma}{\partial t^\gamma}f(t):=\frac{1}{\Gamma(m-\alpha)}\int_0^t \frac{f^{(m)}(u)}{(t-u)^{1+\gamma-m}}du.
\end{equation*}
If $\gamma$ is an integer $\frac{\partial^\gamma}{\partial t^\gamma}f(t)$ is the usual derivative $f^{(\gamma)}$.

\medskip

\noindent {\bf Proof of Theorem \ref{transform}} Remember that $\Gamma(1/2)=\pi$. Then for $\gamma=1/2$,
\begin{equation*}
\frac{\partial^{1/2}}{\partial t^{1/2}}\frac{e^{-\xi^2/4t}}{\sqrt{4\pi t}}=\frac{1}{\sqrt{\pi}}\int_0^t \frac{\partial}{\partial u}\left(\frac{e^{-\xi^2/4u}}{\sqrt{4\pi u}}\right)\cdot \frac{1}{\sqrt{t-u}}du=\frac{\xi e^{-\xi^2/4t}}{4\sqrt{\pi}t^{3/2}}
\end{equation*}
Let $w(t,x):=\E (v(x,|B_{2t}|))=\frac{2}{\sqrt{4\pi t}}\int_0^\infty e^{-\xi^2/4t}v(x,\xi)d\xi$.
Then,
\begin{align*}
\frac{\partial^{1/2}}{\partial t^{1/2}}w(t,x)&=\frac{1}{\sqrt{\pi}}\int_0^t \frac{\partial}{\partial u}\frac{2}{\sqrt{4\pi u}}\int_0^\infty e^{-\xi^2/4u}v(x,\xi)d\xi\cdot \frac{1}{\sqrt{t-u}}du\\
&=\frac{2}{\sqrt{\pi}}\int_0^\infty \int_0^t\frac{\partial}{\partial u}\left(\frac{1}{\sqrt{4\pi u}} e^{-\xi^2/4u}\right)\cdot \frac{1}{\sqrt{t-u}}du\ v(x,\xi)d\xi\\
&=\frac{1}{2\sqrt{\pi}}\int_0^\infty \frac{\xi e^{-\xi^2/4t}}{t^{3/2}} v(x,\xi)d\xi
\end{align*}
From \cite{bragg}, $u(x,t)=\frac{1}{2\sqrt{\pi}t^{3/2}}\int_0^\infty\xi e^{-\xi^2/4t}v(x,\xi)d\xi$
is solution of (\ref{eq:49}).\hfill \framebox[0.6em]{}

\subsection{Time change.}
\begin{thm}\label{changementdetemps}
Let $\alpha:\R \rightarrow {\R}^+$ be a positive increasing and differentiable function. Let $X$ be an $\R^d$ valued process and $Y$ a real valued process satisfying the assumptions of Section \ref{algorithme} such that $\sigma_Y$ and $b_Y$ in $(\ref{edsY})$ depend only on time. Let us define
\begin{equation*}
\tilde{\sigma}_X(t,x):=\sigma_X(t,x)\sqrt{\alpha'(\alpha^{-1}(t))},\quad {\tilde b}_X(t,x):=b_X(t,x)\alpha'(\alpha^{-1}(t)),
\end{equation*}
and $\tilde{a}_X:=\tilde{\sigma}_X \tilde{\sigma}_X^T$. Consider the iterated process $Z:=X(\alpha(Y_t))$. Its density $p_Z(t,x,z)$ satisfies
\begin{equation*}
\frac{\partial}{\partial t}p_Z(t,x,z)=\frac{1}{2}a_Y(t)\Gamma^2 p_Z(t,x,z)+b_Y(t)\Gamma p_Z(t,x,z),
\end{equation*}
where the differential operator $\Gamma$ acts on smooth functions $\varphi$ by
\begin{equation*}
\Gamma\varphi(z):=\frac{1}{2}\sum_{i,j=1}^d\frac{\partial^2}{\partial z_i\partial z_j} \left[\tilde{a}^{i,j}_X(z)\varphi(z)\right]-\sum_{i=1}^d\frac{\partial}{\partial z_i} \left[\tilde{b}^i_X(z)\varphi(z)\right].
\end{equation*}
\end{thm}
For instance, we can take $\alpha(x):=e^x$.\\

\noindent {\bf Proof of Theorem \ref{changementdetemps}} By independence of $X$ and $Y$, the density of $Z$ is given by
$$p_Z(t,x,y)=\int_\R p_Y(t,0,y)p_X(\alpha(y),x,z)dy.$$
We then have
\begin{align*}
\frac{\partial}{\partial t}p_Z(t,x,y)&=\int_\R \frac{\partial}{\partial t}p_Y(t,0,y)p_X(\alpha(y),x,z)dy\\
&=\int_\R \left(a_Y(t)\frac{\partial^2}{\partial y^2}-b_Y(t)\frac{\partial}{\partial y}\right)p_Y(t,0,y)p_X(\alpha(y),x,z)dy\\
&=a_Y(t)\int_\R p_Y(t,0,y)\frac{\partial^2}{\partial y^2}p_X(\alpha(y),x,z)dy\\
&\quad +b_Y(t)\int_\R p_Y(t,0,y)\frac{\partial}{\partial y}p_X(\alpha(y),x,z)dy
\end{align*}
\noindent Now, since $\alpha'(y)a^{i,j}_X(\alpha(y),z)=\tilde{a}^{i,j}_X(z)$ and $\alpha'(y)b^i_X(\alpha(y),z)=\tilde{b}^i_X(z)$, using the forward Kolmogorov (or Fokker-Planck) equation for $p_X$, we are able to conclude since
\begin{align*}
\frac{\partial}{\partial y}p_X(\alpha(y),x,z)&=\alpha'(y)\left.\frac{\partial}{\partial t}p_X(t,x,z)\right|_{t=\alpha(y)}\\
&=\frac{1}{2}\alpha'(y)\sum_{i,j=1}^d\frac{\partial^2}{\partial z_i\partial z_j} \left[a^{i,j}_X(\alpha(y),z)p_X(\alpha(y),x,z)\right]\\
&-\alpha'(y)\sum_{i=1}^d\frac{\partial}{\partial z_i} \left[b^i_X(\alpha(y),z)p_X(\alpha(y),x,z)\right]\\
&=\Gamma p_X(\alpha(y),x,z).
\end{align*}
\hfill \framebox[0.6em]{}

\section{Position process indexed by the real line.}\label{funakipde}

\medskip

\medskip

\noindent The PDEs that we have associated to iterated processes so far exhibit terms depending on the initial value (cf. Theorem \ref{FKitere} where $f(\cdot)\equiv v(0,\cdot)$ and some of its derivatives appear on the right-hand side of (\ref{edpFKitere})). The PDEs obtained in \cite{allouba} and \cite{erkan} for iterated processes also have this drawback due to the use of the absolute value of the time process. Another type of  PDE without this drawback is obtained in \cite{funaki} but the underlying iterated process takes values in the complex plane. Extending the construction (\ref{realline}), we obtain in this section a Feynman-Kac formula where the initial value does not appear any longer in the PDE with a real valued underlying process.

\medskip

\noindent Let us consider $(X_+(t))_{t\geq 0}$ and $(X_-(t))_{t\geq 0}$ two Markov processes with infinitesimal generators ${\cal L}_+$ (resp. ${\cal L}_-$) such that $X_+$ (resp. $X_-$) takes values in $\R^d$ (resp. $\R^n$). Inspired by Funaki's construction \cite{funaki}, we define the $\R^{d+n}$-valued process $(X^{(x_1,x_2)}_t)_{t\in \R}$  starting from $(x_1,x_2)\in \R^{2d}$ and given for any real time index by
\begin{equation}
\label{XindexeparR}
X^{(x_1,x_2)}_t:=\left\lbrace
\begin{array}{lcl} 
(X_+^{x_1}(t),x_2)& {\rm if} \, t\geq 0,\\ 
(x_1,X_-^{x_2}(-t))& {\rm if} \, t<0,
\end{array}\right.
\end{equation}
Note that in \cite{funaki} the resulting process takes values in the complex plane whereas each component of our $X^{(x_1,x_2)}$ is real valued.

\medskip

\noindent We will be interested in bounded functions $f:\R^d\rightarrow \R$ which admit an extension $\tilde{f}:\R^{d+n}\rightarrow \R$ satisfying
\begin{eqnarray}\label{extension}
&&(i)\,\, \,  \forall x_1\in \R^d,\,  \tilde{f}(x_1,0)=f(x_1),\\\
&&(ii)\,\, \,  \forall x_2\in \R^n,\, \tilde{f}(\cdot,x_2)\in {\cal D}({\cal L}_+),\nonumber\\
&&(iii)\,\, \, \forall x_1\in \R^d,\, \tilde{f}(x_1,\cdot)\in {\cal D}({\cal L}_-),\nonumber\\
&&(iv)\,\, \,  \forall (x^1,x^2)\in \R^d\times\R^n, (\mathcal{L_+}+\mathcal{L_-})\tilde{f}(x_1,x_2)=0.\nonumber
\end{eqnarray}
In (iv), the operator $\mathcal{L_+}$ (resp. $\mathcal{L_-}$) acts on $x_1$ (resp. $x_2$). Let us stress the resemblance between (iv) and the intertwining identity (\ref{edppourLambda}) in section \ref{sec:intertwining}.

\bigskip

\subsection{Main statement}
\begin{thm}
\label{edpsansf} 
Let $X$ be defined by (\ref{XindexeparR}). Consider $Y$ a real valued continuous process, independent of $X_+$ and $X_-$, such that $Y_0=0$. Let us assume that $Y$ admits a density $p_Y(t,0,y)$ which satisfies
\begin{equation*}
\frac{\partial}{ \partial t}p_Y(t,0,y)=P(t,\frac{\partial}{ \partial t},\frac{\partial}{ \partial u})p_Y(t,0,y),\, \, \, \forall\,  t>0,\, \forall\,  y\in \R,
\end{equation*}
for some polynomial $P$ with constant coefficients. Suppose moreover that $\frac{\partial^k}{ \partial t^k}p_Y(t,0,\cdot)$ is integrable for all $k$ up to the order of $P$ w.r.t. its second variable. 

Then 
\begin{equation*}\label{sansf}
v(t,x):=\E\left[\tilde{f}(X^{(x,0)}_{Y_t})\right],
\end{equation*}
is solution of the PDE
\begin{equation}\label{FK}
  \frac{\partial}{\partial t}v(t,x)=P(\frac{\partial}{ \partial t},-\, \mathcal{L}_+)v(t,x),\quad \forall \, x\in \R^d,\, \forall\,  t>0,
\end{equation}
with initial condition $v(0,x)=f(x)$, where $\mathcal{L}_+$ acts on $v$ as a function of $x$.
\end{thm}

\bigskip

\noindent It can be noticed that if $X_+$ and $X_-$ are independent and $f$ satisfies assumptions (\ref{extension}), then for all $t>0$ the function $x\mapsto \E \tilde{f}(X_+^x(t),X_-^0(t))$ satisfies (\ref{extension}) too, with an extension given by $(x,y)\mapsto \E \tilde{f}(X_+^x(t),X_-^y(t))$.\\
Theorem \ref{edpsansf} applies if we choose for instance $f(x):=e^{-x^2}$ with its extension $\tilde{f}(x,y):=e^{-x^2}\cdot \arctan(y)$ (cf. (\ref{extension})), and independent diffusions $X_+$ and $X_-$ with infinitesimal generators $\mathcal{L}_+g(x)=\frac{1}{2}\partial_{xx}g(x)+x\partial_xg(x)$ and $\mathcal{L}_-g(y)=(1+y^2)\partial_{yy}+3\partial_y$.

\bigskip

\noindent Let us mention that the assumptions can be weakened and $\tilde{f}$ can be unbounded when $X_-$ or $X_+$ is a diffusion process. For instance, let $X_+$, $X_-$, $Y$ be independent, $X_+$ be an Ornstein-Uhlenbeck process and $X_-$ and $Y$ be two Brownian motions. In this case $\mathcal{L}_+=\frac{1}{2}\frac{\partial^2}{\partial x^2}-x\frac{\partial}{\partial x}$ and $\mathcal{L}_-=\frac{1}{2}\frac{\partial^2}{\partial y^2}$. Let $f(x):=x$ extended in $\tilde{f}(x,y):=x\cosh(y)$. Then $(\mathcal{L}_++\mathcal{L}_-)\tilde f(x,y)=0$ and $\tilde{f}(x,0)=f(x)=x$, $\forall \, x,y\in \R$. Then for $t>0$,
\begin{equation*}
\E\left[\tilde{f}(X_t^{(x,0)})\right]=\E\left[\tilde{f}(X_+^x(t),0)\right]=\E\left[X_+^x(t)\right]=xe^{-t/2}.
\end{equation*}
If $t<0$,
\begin{equation*}
\E\left[\tilde{f}(X_t^{(x,0)})\right]=\E\left[\tilde{f}(x,X_-^0(-t))\right]=\E\left[x\cosh(X_-^0(-t))\right]=xe^{-t/2}.
\end{equation*}
Therefore
\begin{equation*}
\forall\,  t\in \R, \quad \E\left[\tilde{f}(X^{(x,0)}_t)\right]=xe^{-t/2}
\end{equation*}
Then iterating by $Y_t$, and setting $v(t,x):=\E\left[\tilde{f}(X^{(x,0)}(Y_t))\right]]$, we obtain $v(t,x)=\E\left[xe^{-Y_t/2}\right]=xe^{t/8}$ which indeed satisfies 
\begin{eqnarray*}
 \frac{\partial}{\partial t}v(t,x)&=&\frac{1}{2}(\mathcal{L}_+)^2v(t,x)\\
&=&\frac{1}{8}(\partial_x^4-3x\partial_x^3+(x^2-2)\partial_x^2+x\partial_x)v(t,x),
\end{eqnarray*}
forall $x\in \R^d$ and $t>0$, with initial condition $v(0,x)=x$, as stated in (\ref{FK}).

\bigskip

\noindent {\bf Proof  of Theorem \ref{edpsansf}.} For $t\in\R$ and $x\in \R^d$, define $\psi(x,t):=\E\left[\tilde{f}(X^{(x,0)}_t)\right]$. Remember that the notation $X^{(x,0)}$ implies that $X_-$ starts at $0$ this is why we write $X_-^0$ below. We start by proving that 
\begin{equation}\label{FKpourX}
 \frac{\partial}{\partial t}\psi(x,t)=\mathcal{L}_+\psi(x,t),\quad \forall \, t\in \R, \, \, \forall \, x\in\R^d.
\end{equation}
If $t>0$, then $\psi(x,\cdot)=\E(f(X_+^x(\cdot)))$ on an open neighborhood of $t$ and therefore $\frac{\partial}{\partial t}\psi(x,t)=\mathcal{L}_+\psi(x,t).$

\noindent If $t<0$, for all $s$ in an open interval containing $t$, we have $\psi(x,s)=\E(\tilde{f}(x,X_-^0(-s)))=\E(g(X_-^0(-s)))=P_-^{-s}g\, (0)$ with $g: x_2\mapsto \tilde {f}(x,x_2)$ where $x$ plays the role of a parameter. Therefore $\frac{\partial}{\partial t}\psi(x,t)=-\, P_-^{-t}\, {\cal L}_- \, g\, (0).$ The operator ${\cal L}_-$ acts on the second variable $x_2$. With the notations of (\ref{extension}), ${\cal L}_- \, g(x,x_2)$ is also ${\cal L}_- \, \tilde{f}$ so using $(iv)$ of (\ref{extension}), we obtain $-{\cal L}_- \, g(x,x_2)={\cal L}_+\tilde{f}(x,x_2)$ where ${\cal L}_+$ acts only on the first variable $x$. We conclude that $\frac{\partial}{\partial t}\psi(x,t)=P_-^{-t}\, {\cal L}_+\, g\, (0)={\cal L}_+\, P_-^{-t}\, g\, (0)$, the latter identity being true since ${\cal L}_+$ acts only on $x$. 

\noindent It remains to study the case $t=0$. Previously we obtained that $\frac{\partial}{\partial t}\psi(x,t)=\mathcal{L}_+\, P_+^t\,  f\ (x)$ for $t>0\, $ and $\frac{\partial}{\partial t} \psi(x,t)=\mathcal{L}_+\, P_-^{-t}\,  g\, (0)$ for $t<0\, $. The left-hand sides admit the same limit when $t\rightarrow 0$ equal to  $\mathcal{L}_+\, f\, (x)$ since $\mathcal{L}_+$ acts only on $x$. Hence $t\rightarrow \psi(x,t)$ is differentiable at $t=0$ with derivative $\mathcal{L}_+\, f\, (x)$ which coincides with $\mathcal{L}_+\, \psi\, (x,0)$. Hence (\ref{FKpourX}) is proved.

\noindent We now consider $v(t,x):=\E\left[\tilde{f}(X^{(x,0)}_{Y_t}\right]$. Then $v(t,x)=\int_\R p_Y(t,0,s)\psi(x,s)ds$ by independence of $X$ and $Y$ and the following identities hold,
\begin{eqnarray*}
\frac{\partial}{\partial t}v(t,x)&=&\frac{\partial}{\partial t}\int_\R p_Y(t,0,s)\psi(x,s)ds\\
&=&\int_\R P(t,\frac{\partial}{\partial t},\frac{\partial}{\partial s})p_Y(t,0,s)\psi(x,s)ds.
\end{eqnarray*}

If we perform successive integrations by parts on each term involving $\frac{\partial}{\partial s}$ and its powers that appear in $P$, we see that ${(-1)}^k\, \int_\R \frac{\partial^k}{\partial s^k}p_Y(t,0,s)\psi(x,s)ds=\int_\R p_Y(t,0,s)\frac{\partial^k}{\partial s^k}\psi(x,s)ds$ for all integer $k$. Then we conclude using (\ref{FKpourX}) that
\begin{equation*}
\int_\R \frac{\partial^k}{\partial s^k}p_Y(t,0,s)\psi(x,s)ds=(-1)^k\, \int_\R p_Y(t,0,s)({\cal L_+})^k\psi(x,s)ds
\end{equation*} 
where ${\cal L_+}$ acts on $x$. In this way we have been able to separate the variables $t$ and $x$. Therefore
\begin{equation*}
\frac{\partial}{\partial t}\E\left[\tilde{f}(X^{(x,0)}(Y_t)\right]=P(t,\frac{\partial}{\partial t},-\, \mathcal{L}_+) \int_\R p_Y(t,0,s)\psi(x,s)ds.
\end{equation*}
This is  (\ref{FK}) that we wanted to prove. \hfill \framebox[0.6em]{}

\subsection{Application of Theorem \ref{edpsansf}}

\begin{corol}\label{FKwithkilling}

Let us keep the assumptions and notations of Theorem \ref{edpsansf}. Furthermore, we assume $X_-$ is such that the set $\{s\in [0,t], X_-(s)=0\}$ has Lebesgue measure zero $\P$-almost surely for all $t>0$. Let $c_+:\R^d\rightarrow\R_+$ and $c_-:\R^n\rightarrow\R_+$ be two continuous bounded functions. Let $f$ satisfying assumptions (i) to (iii) of (\ref{extension}) and the relation 
$$\forall (x^1,x^2)\in \R^d\times\R^n, (\mathcal{L_+}-c_+(x_1)+\mathcal{L_-}-c_-(x_2))\tilde{f}(x_1,x_2)=0.$$
Set $c(x_1,x_2)=c_+(x_1)\1_{x_2=0}-c_-(x_2)\1_{x_2\neq 0}$, then \\$\nu(t,x):=\E \left[ \exp \left\{ -\int_0^{Y_t} c(X^{(x,0)}_s)ds\right\}\tilde{f}(X^{(x,0)}_{Y_t})\right]$ is solution of the PDE
\begin{equation}\label{FKkilling}
\frac{\partial}{\partial t}\nu(t,x)=P(\frac{\partial}{\partial t},-(\mathcal{L}_+-c_+(x)))\nu(t,x), \, \, \, \forall \, x\in \R^d,\, \forall\,  t>0,
\end{equation}
with initial condition $\nu(0,x)=f(x)$.

\end{corol}

\bigskip

\noindent In the following statement we show that the Euler-Bernoulli {\it beam equation}
\begin{equation}
\label{eq:47}
\frac{\partial^2}{\partial x^2}\left(g(x)\frac{\partial^2u}{\partial x^2}\right)+m(x)\frac{\partial^2 u}{\partial t^2}=0, \quad t>0,\, \,  0<x<L,
\end{equation}
where $g(x)>0$ is the flexural rigidity and $m(x)>0$ the lineal mass density, can be obtained as a consequence of Theorem \ref{edpsansf}, by considering the iteration of a Brownian motion by an independent Cauchy process.

\begin{corol}\label{beam}
Let $g$ and $m$ be two positive functions. Let $X_+$ be a diffusion process with infinitiesmal generator $\mathcal{L}_+:=g(x)\partial_{xx}^2$. Consider $X_-$ (resp. C) a Markov (resp. Cauchy) process such that $X_+$, $X_-$ and $C$ are independent. Define $\gamma_t:=C_{t/\sqrt{g(x)m(x)}}$. Then for any $f$ extendable in $\tilde{f}$ in the sense of (\ref{extension}), the function 
\begin{equation*}
u(t,x):=\E[\tilde{f}(X^{(x,0)}_{\gamma_t})],
\end{equation*}
satisfies equation $(\ref{eq:47})$ with initial condition $u(0,x)=f(x), \forall \, x\in \R$. 
\end{corol}

\medskip

\noindent As an example, let be $Y$ a Brownian motion with drift $\mu$. Corollary \ref{FKwithkilling} provides a probabilistic representation of the solution to equation
$$
\left \lbrace 
\begin{array}{lcl} 
\frac{\partial}{\partial t}u(t,x)=\frac{1}{2}\left(\mathcal{L}_+-c_+\right)^2+\mu \left(\mathcal{L}_+-c_+\right), \quad \forall x\in \R,\forall t>0,\\
u(0,x)=f(x),
\end{array}\right.
$$
by $u(t,x)=\E \left[ \exp \left\{ \int_0^{Y_t} c(X^{(x,0)}(s))ds\right\}\tilde{f}(X^{(x,0)}_{Y_t})\right]$. A solution to such an equation can be processed using the algorithm developed in this paper and the remark following Theorem \ref{rateofconvergence}.

\bigskip

\noindent {\bf Proof of Corollary \ref{FKwithkilling}:} Let us define $(\tilde{X}_+(t))_{t\geq 0}$ (resp. $(\tilde{X}_-(t))_{t\geq 0}$) two Markov processes with infinitesimal generators ${\cal L}_+-c_+$ (resp. ${\cal L}_--c_-$) and the process
\begin{equation}
\tilde{X}^{(x_1,x_2)}_t:=\left\lbrace
\begin{array}{lcl} 
(\tilde{X}_+^{x_1}(t),x_2)& {\rm if} \, t\geq 0,\\ 
(x_1,\tilde{X}_-^{x_2}(-t))& {\rm if} \, t<0,
\end{array}\right.
\end{equation}
From Theorem \ref{edpsansf}, $\nu(t,x):=\E\left[\tilde{f}(\tilde{X}^{(x,0)}_{Y_t})\right]$ is solution of 
$$\frac{\partial}{\partial t}\nu(t,x)=P(\frac{\partial}{ \partial t},-\, (\mathcal{L}_+-c_+(x)))\nu(t,x)$$
For $t\geq 0$, using the definition (\ref{XindexeparR}) of $X^{(x_1,x_2)}$ and the Feynman-Kac formula,
\begin{multline*}
\E\left[\tilde{f}(\tilde{X}^{(x,0)}_t)\right]=\E\left[\tilde{f}(\tilde{X}_+^x(t),0)\right]\\=\E \left[\exp\left\{-\int_0^tc_+(X_+^x(s))ds\right\}\tilde{f}(X_+^x(t),0)\right]\\
=\E \left[\exp\left\{-\int_0^tc(X^{(x,0)}(s))ds\right\}\tilde{f}(X^{(x,0)}_t)\right]\
\end{multline*}
When $t<0$, the following identities hold
\begin{align*}
\E\left[\tilde{f}(\tilde{X}^{(x,0)}_t)\right]&=\E\left[\tilde{f}(x,\tilde{X}_-^0(-t))\right]\\
&=\E \left[\exp\left\{-\int_0^{-t}c_-(X_-^0(s))ds\right\}\tilde{f}(x,X_-^0(-t))\right]\\
&=\E \left[\exp\left\{\int_0^t c_-(X_-^0(-s))ds\right\}\tilde{f}(x,X_-^0(-t))\right]\\
\end{align*}
Since $c(X^{(x,0)}(s)=-c_-(X_-^0(-s))$ if $X_-^0(-s)\neq 0$, using assumption on the zero set of $X_-$, we have
$$\E\left[\tilde{f}(\tilde{X}^{(x,0)}_t)\right]=\E \left[\exp\left\{-\int_0^tc(X^{(x,0)}(s))ds\right\}\tilde{f}(X^{(x,0)}_t)\right]$$
So 
$$\nu(t,x)=\E\left[\tilde{f}(\tilde{X}^{(x,0)}_{Y_t})\right]=\E \left[\exp\left\{-\int_0^{Y_t}c(X^{(x,0)}(s))ds\right\}\tilde{f}(X^{(x,0)}_{Y_t})\right]$$
\hfill \framebox[0.6em]{}

\medskip

\noindent {\bf Proof of Corollary \ref{beam}.} Since the density of $C$, $p_C(t,0,y):=\P(C_t\in dy)=\frac{t}{\pi(t^2+y^2)}$ satisfies 
$$\frac{\partial^2}{\partial t^2}p_C(t,0,y)=-\frac{\partial^2}{\partial y^2}p_C(t,0,y),$$
$w(t,y):=p_C(t/\sqrt{g(x)m(x)},0,y)$ is solution of
$$g(x)m(x)\frac{\partial^2}{\partial t^2}w(t,y)=-\frac{\partial^2}{\partial y^2}w(t,y).$$
Applying Theorem \ref{edpsansf} gives $w$ is solution of
$$(\mathcal{L}_+)^2+g(x)m(x)\frac{\partial^2 u}{\partial t^2}=0.$$
Factorizing by $g(x)$ ends the proof. \hfill \framebox[0.6em]{}

\section{Application of the numerical scheme to the iterated Brownian motion}\label{application}

\noindent In this section, we illustrate the algorithm proposed in Section \ref{algorithme}. First we simulate a trajectory of an iterated Brownian motion (IBM) $Z_t=X(|Y_t|)$ on $[0,T]$ where $X$ and $Y$ are two independent Brownian motions. Then we approximate numerically the function $v(t,x):=\E \left[f(Z^x_t))\right]$  and the variations of order three and four of $(Z_t)$.

\smallskip

\begin{figure}[h]
\begin{center}
\includegraphics[trim = 10mm 9mm 0mm 9mm, clip, scale=0.7]{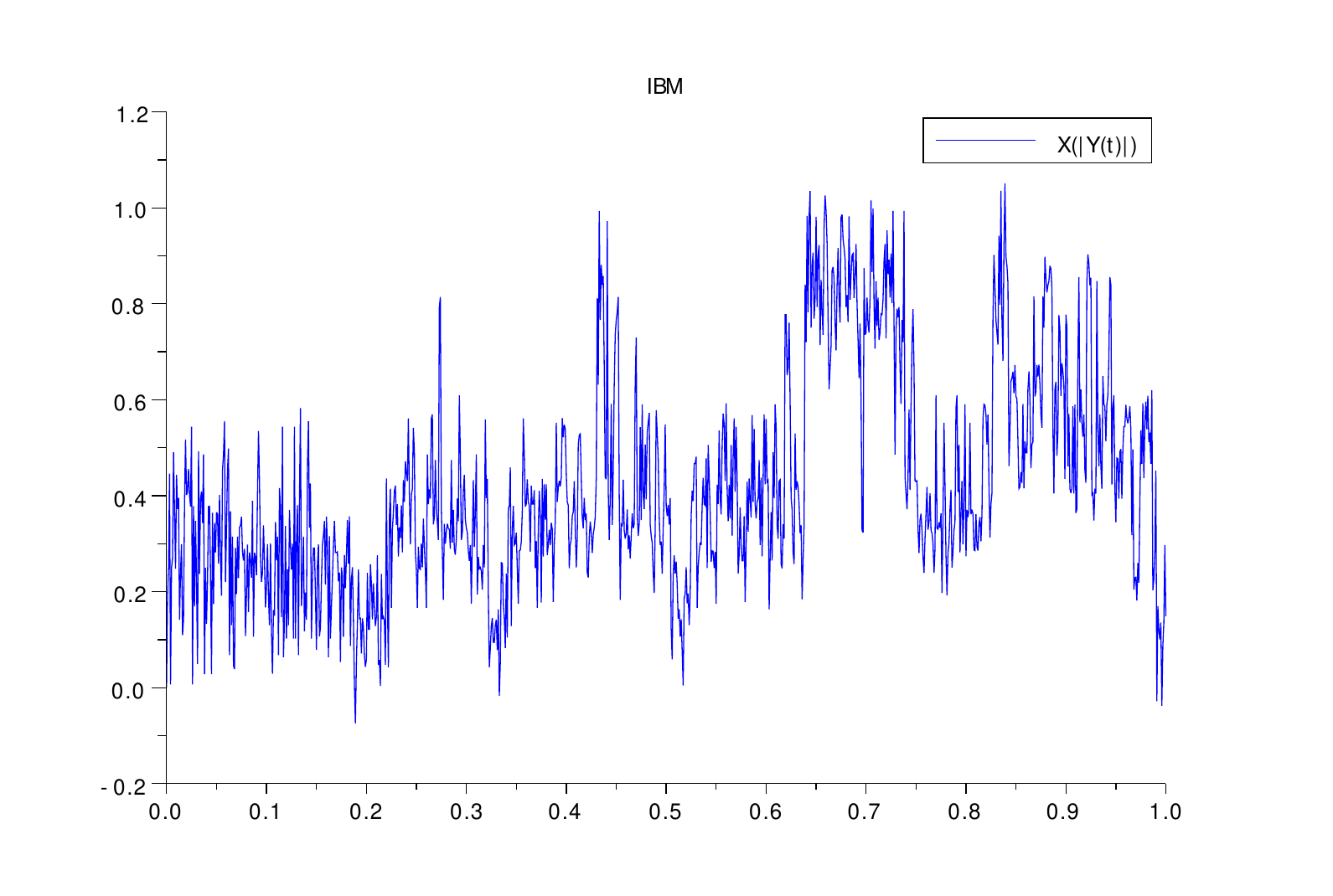}
\end{center}
\caption{\label{fig:1} A trajectory of $\tilde{Z}^n$ on $[0,1]$ for $n=1000$.}
\begin{center}
\includegraphics[trim = 18mm 50mm 20mm 50mm, clip,scale=0.7]{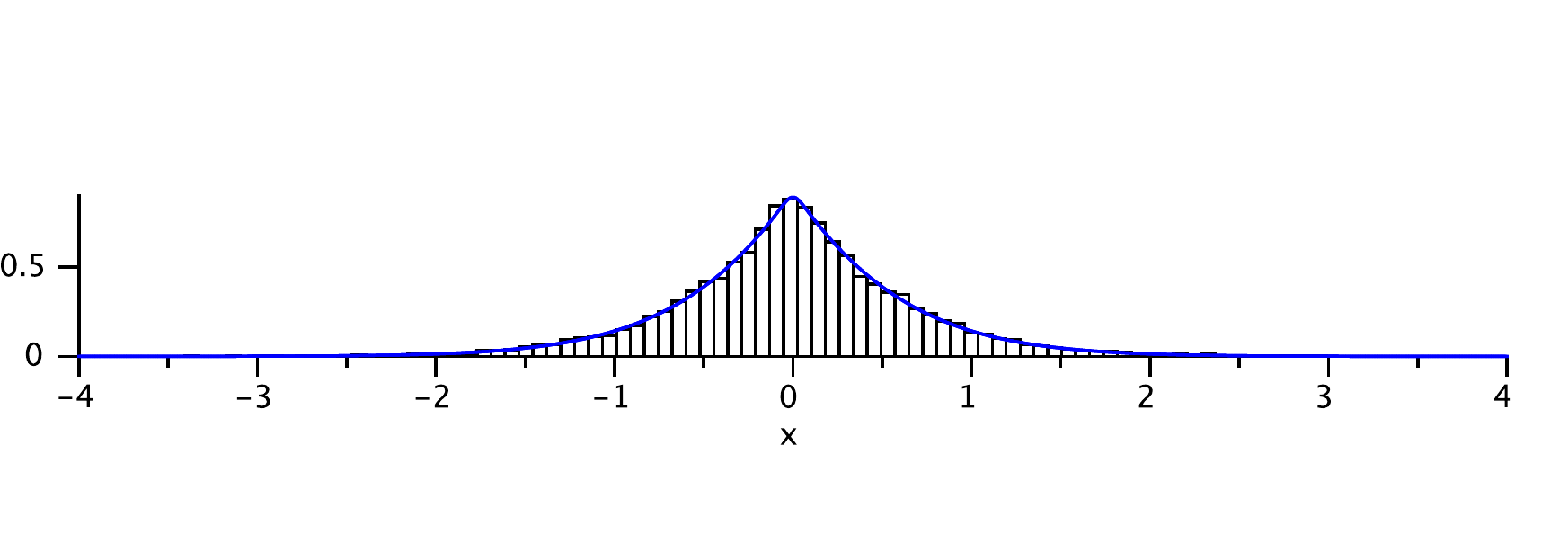}
\end{center}
\caption{\label{fig:2}Comparison between the density of the Iterated Brownian Motion (IBM) at $t=1$ ($Z_1$) and an histogram representing an estimation of this density by the simulation of $20000$ trajectories of $\tilde{Z}_1$ for $n=1000$.}
\end{figure}

\begin{figure}[h]
\begin{center}
\includegraphics[trim = 10mm 9mm 0mm 9mm, clip, scale=0.7]{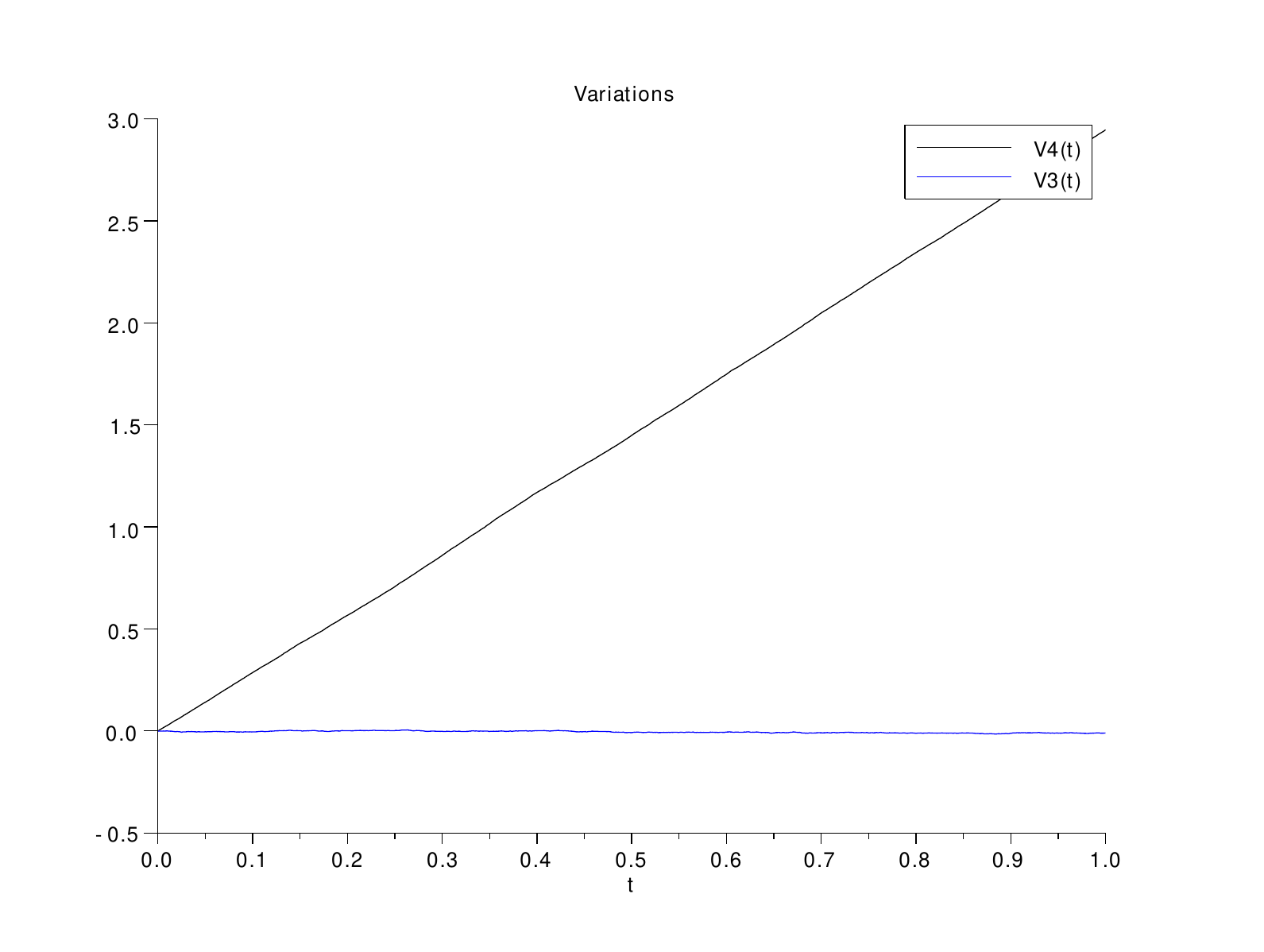}
\end{center}
\caption{\label{fig:3} Estimation of $V_3$ and $V_4$ from the simulation of $2000$ trajectories of $\tilde{Z}^n$ for $n=1000$.}
\end{figure}

\begin{figure}[h]
\begin{center}
\includegraphics[trim = 1mm 9mm 0mm 9mm, clip, scale=0.6]{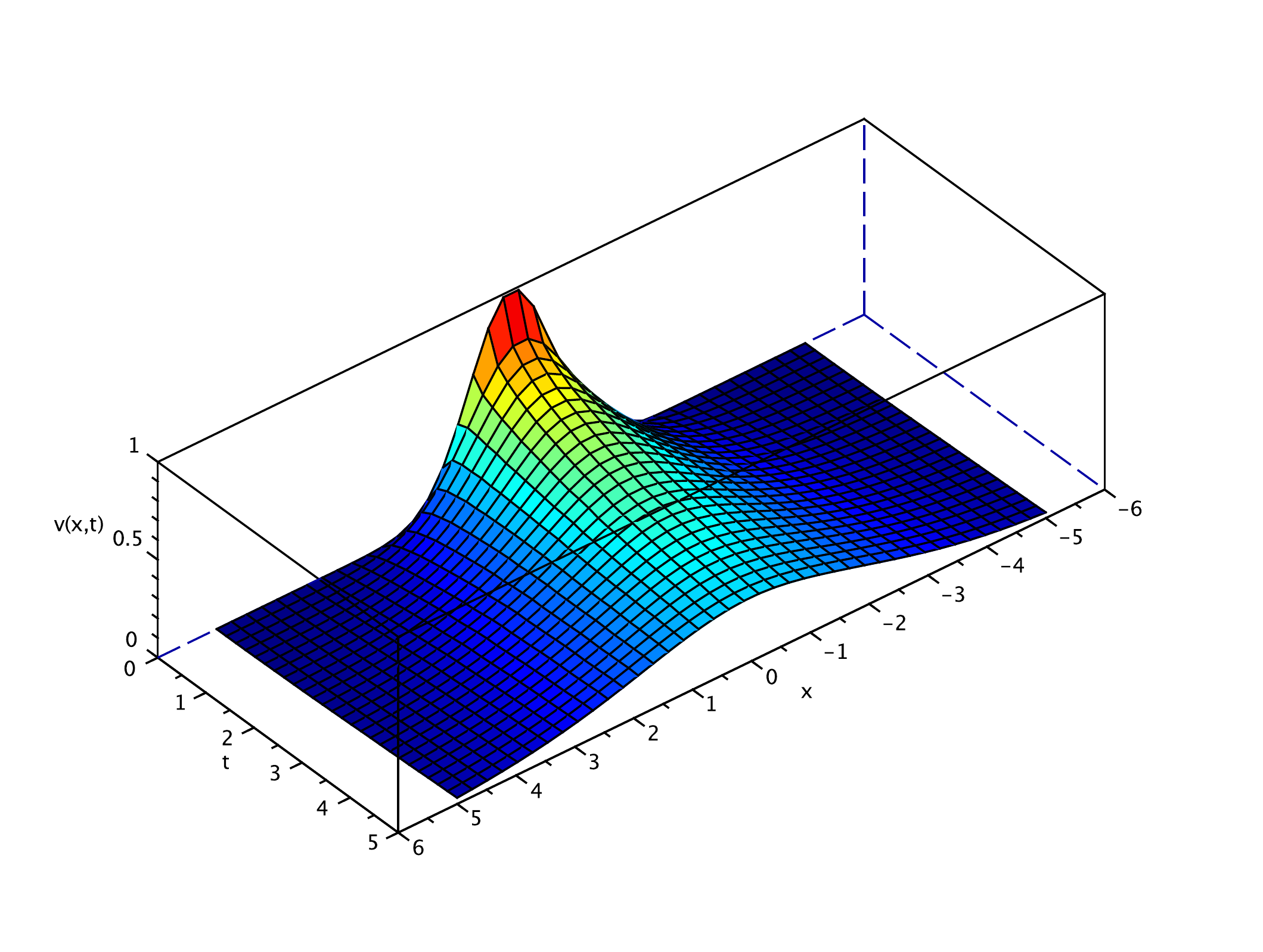}
\end{center}
\caption{\label{fig:5} Estimation of $v(t,x)=\E \exp\{-(X^x(|Y(t)|))^2\}$, solution of (\ref{eq:45}). ($n=1000$, $20000$ trajectories)} 
\end{figure}

\noindent For a fixed positive integer $n$, the Brownian motion $Y$ is evaluated at times $kT/n$ and we define our piecewise constant processes $(\bar{Y}_t^n)_{0\leq t \leq T}$ recursively by
\begin{equation*}
\bar{Y}_{(k+1)T/n}^n=\bar{Y}_{kT/n}^n+\xi_k^n,
\end{equation*}
where $\xi_k^n$ are independent centered Gaussian random variables with variance $T/n$. This determines $M_n=\sup_{t\in [0,T]} |\bar{Y}_t^n|$ which is a.s. finite. The same construction for $\bar{X}_t^n$ on $[0,M_n]$ can be performed, 
$$\bar{X}_{(k+1)M_n/n}^n=\bar{X}_{kM_n/n}^n+\zeta_k^n$$
where $(\zeta_k^n)$ are independent centered Gaussian random variables with variance $M_n/n$.\\
The composition $\bar{X}^n(|\bar{Y}^n(kT/n)|)$ for $k=0,1,\ldots, n$, is given by
\begin{equation}\label{leschema}
\bar{X}^n(|\bar{Y}^n(kT/n|)=\bar{X}^n\left( \frac{M_n}{n} \Big \lfloor \frac{n}{M_n}|\bar{Y}^n(kT/n)|\Big \rfloor \right),
\end{equation}
with $\lfloor \cdot \rfloor$ the floor function. The continuous approximation $\tilde{Z}^n$ of $Z$ is the linear interpolation between the points defined in (\ref{leschema}). A trajectory of $\tilde{Z}^n$  can be seen in Figure \ref{fig:1} presenting huge variations.\\
As mentioned in section \ref{algorithme} the scheme (\ref{leschema}) is interesting since it converges uniformly as $n$ tends to infinity, it simply requires to simulate $2n$ independent Gaussian random variables and the composition of $\bar{X}^n$ by $|\bar{Y}^n|$ is facilitated by the use of step functions. This makes it possible to use methods such as Monte Carlo methods which require the simulation of thousands of trajectories (see Figure \ref{fig:2}). Moreover the almost sure uniform convergence proved in this paper for $Z^n$ and its version $\tilde{Z}^n,$  which is more convenient for implementation, renders possible all types of numerical studies requiring an approximation of a whole trajectory and not only of the value at a fixed time, like for instance the variations of various order. In Figure \ref{fig:3} we apply this remark to the third and fourth order variations of $(Z_t)$, illustrating the following results of \cite{burdzy3}
\begin{equation*}
V_3(t):=\lim_{|\Lambda|\rightarrow 0} \sum_{k=1}^n(Z(t_k)-Z(t_{k-1}))^3= 0 \quad \text{in\ } L^p.
\end{equation*}
\begin{equation*}
V_4(t):=\lim_{|\Lambda|\rightarrow 0} \sum_{k=1}^n(Z(t_k)-Z(t_{k-1}))^4= 3t \quad \text{in\ } L^p,
\end{equation*}
where $\Lambda$ denotes an arbitrary subdivision $t_0=0\leq t_1 \leq \ldots \leq t_n=t $ of $[0,t]$ with mesh $|\Lambda|:= \max_{1\leq k \leq n} |t_k-t_{k-1}|$:

\bigskip

\noindent The algorithm can also be used to simulate the solution to a fourth-order PDE of the type we studied in the previous sections. Figure \ref{fig:5} shows an approximation of $v(t,x):=\E \left[f(Z^x_t))\right]$ corresponding to $f(x)=e^{-x^2}$. From Corollary \ref{corol:2} we know that $v(t,x)$ is the unique solution of
\begin{equation}
\label{eq:45}
\frac{\partial}{\partial t}v(t,x)=\frac{1}{2\sqrt{2\pi t}}\frac{\partial^2}{\partial x^2}f(x)+\frac{1}{8}\frac{\partial^4}{\partial x^4}v(t,x), \quad t>0,\ x\in\R^d,
\end{equation}
with initial condition $v(0,x)=f(x)$. 

\bigskip

\section{Appendix}\label{appendix}
\subsection{Classical results for section \ref{algorithme}}
\begin{prop}(cf. \cite{friedman}).
\label{estimeesclassiques} 
Let $T>0$ and $(U_t)_{0\leq t\leq T}$ be the solution of (\ref{edsX}) under assumptions (\ref{lipschitz_cond})-(\ref{holder_cond}). Then for all $t, s\in [0,T]$ and $p\geq 1$ such that $\E|U_0|^{2p}<\infty$,
\begin{equation}
\label{eq:11}
\E(|U_t|^{2p})\leq \left(1+\E(|U_0|^{2p})\right)e^{Ct},
\end{equation}
\begin{equation}
\label{eq:17}
\E(\sup_{t\in [0,T]}|U_t|^{2p})\leq C(1+T^p)\, [\E(|U_0|^{2p})+\left(1+\E(|U_0|^{2p})\right)T^p e^{CT}],
\end{equation}
\begin{equation}
\label{Kolmogorovcriterion}
\E(|U_t-U_s|^{2p})\leq C\left(1+\E(|U_0|^{2p})\right)(1+T^p)|t-s|^p e^{CT},
\end{equation}
for some constant $C>0$ depending only on $K$ and $p$.
\end{prop}

\bigskip

\noindent\textbf{Proof of Proposition \ref{errorboundforX}.}
$C$ can denote different constants from lines to lines depending only on $K$, $p$, $n$ and $\E|X_0|^{2p}$.
Keeping the notations in \cite{faure}, we have for $t\in [t_k, t_{k+1}]$:
\begin{multline}
\label{eq:14}
\E(|\epsilon_t^n|^{2p})\leq \E(|\epsilon_{t_k}^n|^{2p})+C\int_{t_k}^t \E\{ |\epsilon_s^n|^{2p}+|\sigma(s,X_s)-\sigma(t_k,\tilde{X}_{t_k}^n)|^{2p}\\
+|b(s,X_s)-b(t_k,\tilde{X}_{t_k}^n)|^{2p}\}ds
\end{multline}
where $\epsilon_{t_k}^n$ denotes the error process defined by $\epsilon_t^n:=X_t-X_t^n$.
\begin{eqnarray}
\label{eq:13}
|b(s,X_s)-b(t_k,\tilde{X}_{t_k}^n)|^{2p}&\leq& C|b(s,X_s)-b(t_k,X_s)|^{2p}+C|b(t_k,X_s)-b(t_k,X_{t_k})|^{2p}\nonumber\\
&+&C|b(t_k,X_{t_k})-b(t_k,\tilde{X}_{t_k}^n)|^{2p}.
\end{eqnarray}
From Proposition \ref{estimeesclassiques} and the assumptions of Section \ref{algorithme}, we have
\begin{align*}
&\E|b(s,X_s)-b(t_k,X_s)|^{2p}\leq K^{2p}(\Delta_t)^{2p\beta},\\
&\E|b(t_k,X_s)-b(t_k,X_{t_k})|\leq K^{2p}\E|X_s-X_{t_k}|^{2p}\leq C(1+T^p)(\Delta_t)^pe^{CT},\\
&\E|b(t_k,X_{t_k})-b(t_k,\tilde{X}_{t_k}^n)|^{2p}\leq K^{2p}\E|\epsilon_{t_k}^n|^{2p}.
\end{align*}
Combining these inequalities with (\ref{eq:13}) yields 
\begin{equation}
\label{eq:15}
\E|b(s,X_s)-b(t_k,\tilde{X}_{t_k}^n)|^{2p}\leq C\{(1+T^p)(\Delta_t)^\gamma e^{CT}+\E|\epsilon_{t_k}^n|^{2p}\}
\end{equation}
for $\gamma=\max(p,2p\beta)$ if $\Delta_t\geq1$ and $\gamma=\min(p,2p\beta)$ otherwise. The same inequality holds for $\sigma$. By writing $\epsilon_s^n=\epsilon_{t_k}^n+(X_s-X_{t_k})+(\tilde{X}_{t_k}^n-\tilde{X}_s^n)$, and using inequality (\ref{eq:11}), we obtain
$$\E|\epsilon_s^n|^{2p}\leq C\E|\epsilon_{t_k}^n|^{2p}+C(1+T^p)(\Delta_t)^pe^{CT}.$$
Thus, inequality (\ref{eq:14}) for $t=t_{k+1}$ becomes
$$\E(|\epsilon_{t_{k+1}}^n|^{2p})\leq \E(|\epsilon_{t_k}^n|^{2p})(1+C\Delta_t)+C(1+T^p)(\Delta_t)^{\gamma+1} e^{CT}.$$
Consequently,
\begin{align}
\E(|\epsilon_{t_k}^n|^{2p})&\leq Cne^{nC\Delta_t}(1+T^p)(\Delta_t)^{\gamma+1} e^{CT}\notag\\
&\leq C(1+T^p)T(\Delta_t)^\gamma e^{CT}. \label{eq:16}
\end{align}
Using
\begin{multline*}
\epsilon_t^n=\int_0^t \sum_{k=0}^{n-1}(\sigma(s,X_s)-\sigma(t_k,\tilde{X}_{t_k}^n))\1_{[t_k,t_{k+1}]}dW_s\\
+\int_0^t \sum_{k=0}^{n-1}(b(s,X_s)-b(t_k,\tilde{X}_{t_k}^n))\1_{[t_k,t_{k+1}]}(s)ds,
\end{multline*}
and BurkholderDavisGundy inequality, we have
\begin{multline*}
\E (\sup_{t\in [0,T]} |\epsilon_t^n|^{2p})\leq C \E(\int_0^T \sum_{k=0}^{n-1}|\sigma(s,X_s)-\sigma(t_k,\tilde{X}_{t_k}^n)|^{2p}\1_{[t_k,t_{k+1}]}(s)ds)\\
+C\E(\int_0^T \sum_{k=0}^{n-1}|b(s,X_s)-b(t_k,\tilde{X}_{t_k}^n)|^{2p}\1_{[t_k,t_{k+1}]}(s)ds).
\end{multline*}
Using this latter inequality as well as (\ref{eq:15}) and (\ref{eq:16}) we complete the proof.
\hfill \framebox[0.6em]{}\\

\noindent\textbf{Proof of Lemma \ref{lem:1}.}
We suppose $\forall \alpha<l, n^\alpha \sup_{t\in[0,T]} |f_n(t)-f(t)|\rightarrow_n 0$. The second implication can be shown similarly. Let $\alpha>0$ and $t\in [0,T]$, then
\begin{equation}
\label{eq:2}
\sup_{t\in[0,T]}|\bar{f}_n(t)-f(t)|\leq \sup_{t\in[0,T]}|\bar{f}_n(t)-f_n(t)|+\sup_{t\in[0,T]}|f_n(t)-f(t)|.
\end{equation}
By construction,
$$\sup_{t\in[0,T]}|\bar{f}_n(t)-f_n(t)|=\max_k \sup_{t\in [kT/n,(k+1)T/n]}|f_n(t)-f_n(kT/n)|.$$
Let $t, t'\in[0,T]$,
\begin{align}
|f_n(t)-f_n(t')|&\leq |f_n(t)-f(t)|+|f(t)-f(t')|+|f(t')-f_n(t')|\notag \\
&\leq 2\sup_{t\in[0,T]}|f_n(t)-f(t)|+|f(t)-f(t')|.\label{eq:4}
\end{align}
$f$ being $\beta$-H\" older continue for all $\beta<l$, we choose $\beta$ such that $l>\beta>\alpha$ and 
inequality (\ref{eq:4}) becomes for some $C>0$
\begin{equation}
\label{eq:5}
|f_n(t)-f_n(t')|\leq 2\sup_{t\in[0,T]}|f_n(t)-f(t)|+C|t-t'|^\beta.
\end{equation}
Thus,
$$\sup_{t\in [kT/n,(k+1)T/n]}|f_n(t)-f_n(kT/n)|\leq 2\sup_{t\in[0,T]}|f_n(t)-f(t)|+C(T/n)^\beta.$$
Consequently, the right side of inequality (\ref{eq:2}) multiplied by $n^\alpha$ tends to $0$ as $n$ tends to infinity and the result follows.\hfill \framebox[0.6em]{}

\subsection{Classical results for Theorem \ref{FKitere}}

\noindent\textbf{Proof of Lemma \ref{IPP}.}
We first assume that $g$ is infinitely differentiable with compact support. Let $h>0$,
\begin{align*}
\frac{P^h-I}{h}F(g)&=\frac{1}{h}\int_0^\infty g(s)(P^{s+h}-P^s)f(x)ds\\
&=\frac{1}{h}\int_h^\infty (g(s-h)-g(s))P^sf(x)ds-\frac{1}{h}\int_0^hg(s)P^sf(x)ds
\end{align*}
We have to show the convergence of the right hand side as $h\downarrow 0$ in \\$(B(\R,\R), {||\cdot||_\infty})$. Since $P^tf$ converges to $f$ uniformly as $t\downarrow 0$ and $g$ is continuous, $\exists \delta>0$ such that for all $s\in [0,\delta]$, $\| P^sf-f\|_\infty<\epsilon/\|g\|_\infty$ and $|g(s)-g(0)|<\epsilon/\|f\|_\infty$ so
\begin{multline*}
|g(s)P^sf(x)-g(0)f(x)|\leq |g(s)P^sf(x)-g(s)f(x)|+|g(s)f(x)-g(0)f(x)|\\
\leq \|g\|_\infty |P^sf(x)-f(x)|+\| f\|_\infty |g(s)-g(0)|<2\epsilon
\end{multline*}
Thus, $\forall \epsilon>0$, $\exists \delta$ such that $\forall\ 0<h<\delta$,
\begin{align*}
\left\| \frac{1}{h}\int_0^h g(s)P^sf-g(0)fds \right\|_\infty \leq \frac{1}{h}\int_0^h \|g(s)P^sf-g(0)f\|_\infty ds<2\epsilon
\end{align*}
By dominated convergence theorem, using the bound $\|P^sf\|_\infty\leq \|f\|_\infty$ and the regularity of $g$, we have
$$\lim_{h\downarrow 0} \frac{P^h-I}{h}F(g)=-\int_0^\infty g'(s)P^sfds-g(0)f=-F(g')-g(0)f$$
This ensures that $F(g)\in D(\mathcal{L})$ and $\mathcal{L}F(g)=-F(g')-g(0)f$. Let $(\alpha_n)_{n\in \N}$ be an infinitely differentiable mollifier. We define $g_m(x):=g(x)\1_{|x|\leq m}$ and $\tilde{g}_m(x):=g'(x)\1_{|x|\leq m}$. Then for all $n,m\in \N, (g_m \star \alpha_n)$ is $\mathcal{C}^\infty$ with compact support. Since $g_m\in L^1(\R)$ and $g'_m\in L^1$, $(g_m\star \alpha_n)\rightarrow g_m$ in $L^1(\R)$ and $(g_m\star \alpha_n)'=(\tilde{g}_m\star \alpha_n)\rightarrow \tilde{g}_m$ in $L^1(\R)$
$g$ being continuous, $g_m\in L^\infty(\R)$ for all $m\in \N$, so $|(g_m\star \alpha_n)(0)-g(0)|=|(g_m\star \alpha_n)(0)-g_m(0)|\rightarrow_n 0$. This and the convergence of $g_m$ to $g$ and $\tilde{g}_m$ to $g'$ in $L^1(\R)$ imply that we can construct a sequence $(\phi_n)_n$ of $\mathcal{C}^\infty$ functions with compact support such that
$$\phi_n\rightarrow_n g \text{\ in\ } L^1(\R),\quad  \phi'_n\rightarrow_n g'  \text{\ in\ } L^1(\R) \text{\quad and\quad}|\phi_n(0)-g(0)| \rightarrow_n 0$$
For all $n\in \N$, $F(\phi_n)\in D(\mathcal{L})$ and $\mathcal{L}F(\phi_n)=-F(\phi_n')-\phi_n(0)f$.
From the inequality
\begin{eqnarray*}
|F(\phi_n)(x)-F(g)(x)|&=&\left| \int_0^\infty (\phi_n(s)-g(s))P^sf(x)ds\right|\\
&\leq& \|f\|_\infty \int_0^\infty |\phi_n(s)-g(s)|ds,
\end{eqnarray*}
we have $\|F(\phi_n)-F(g)\|_\infty\rightarrow_n 0$ and similarly $\|F(\phi'_n)-F(g')\|_\infty\rightarrow_n 0$.
$$\|\phi_n(0)f(x)-g(0)f(x)\|_\infty \leq \|f\|_\infty |\phi_n(0)-g(0)|\rightarrow_n 0$$

\noindent We have shown that $(F(\phi_n))$ is then a sequence of $D(\mathcal{L})$ converging to $F(g)$ and $\mathcal{L}F(\phi_n)=-F(\phi'_n)-\phi_n(0)f$ is converging to $-F(g')-g(0)f$. Since $\mathcal{L}$ is a closed operator, we conclude that $F(g)\in D(\mathcal{L})$ and $\mathcal{L}F(g)=-F(g')-g(0)f$.
The second part of the Lemma is proved by induction. \hfill \framebox[0.6em]

\begin{lem}\label{intervertion}
Let $g$ be a continuous density function. For all $f$ in the domain of $\mathcal{L}$,
$$\mathcal{L}\int_0^\infty P^sf(x)g(s)ds=\int_0^\infty \mathcal{L}P^sf(x)g(s)ds.$$
\end{lem}

\noindent\textbf{Proof of Lemma \ref{intervertion}.}
Let $h>0$ and $p(h,x,dy):=\P(X_h^x \in dy)$
\begin{align*}
P^h\int_0^\infty P^sf(x)g(s)ds&=\int_\R p(h,x,dy)\int_0^\infty P^sf(y)g(s)ds\\
&=\int_0^\infty \int_\R p(h,x,dy) P^sf(y)g(s)ds\\
&=\int_0^\infty P^{h+s}f(x)g(s)ds
\end{align*}
So
\begin{equation}
\label{eq:39}
\frac{P^h-I}{h}\int_0^\infty P^sf(x)g(s)ds=\int_0^\infty \frac{P^{h+s}-P^s}{h}f(x)g(s)ds
\end{equation}
where $I$ stands for the identity.
Since, $(P^{h+s}f(x)-P^sf(x))/h\rightarrow \mathcal{L}P^sf(x)$ uniformly in $x$ as $h\downarrow 0$, letting $h$ decreasing to $0$ in equality (\ref{eq:39}) ends the proof.

\hfill \framebox[0.6em]

\nocite{*}
\newpage
\renewcommand{\refname}{Bibliography}
\bibliographystyle{plain}
\bibliography{bibliographie}

\end{document}